\documentclass[reqno]{amsart}
\usepackage[utf8]{inputenc}
\usepackage{mathtools}
\usepackage{dsfont}
\usepackage{a4wide,amsmath}
\usepackage{mathrsfs}
\usepackage{amssymb}
\usepackage{amsthm}
\usepackage{caption}
\usepackage{subcaption}
\usepackage{booktabs}
\usepackage{enumerate}
\numberwithin{equation}{section}
\numberwithin{figure}{section}
\numberwithin{table}{section}
\makeatletter \let\c@table\c@figure \makeatother %use same counter for tables and figures, for less confusing references to them
\usepackage[section]{placeins}
\usepackage{hyperref}
\usepackage[dvipsnames]{xcolor}
\hypersetup{hidelinks}

\newcommand{\R}{\mathbb{R}}
\newcommand{\N}{\mathbb{N}}
\newcommand{\C}{\mathbb{C}}
\newcommand{\Z}{\mathbb{Z}}
\newcommand{\I}{\mathrm{i}}    % imaginary unit
\newcommand{\e}{\mathrm{e}}    % Euler's number
\newcommand{\di}{\mathrm{d}}   % differential
\newcommand{\redel}{\mathop{\textup{Re}}}
\newcommand{\imdel}{\mathop{\textup{Im}}}
\newcommand{\supp}{\mathop{\textup{supp}}}
\newcommand{\norm}[1]{\lVert#1\rVert}
\newcommand{\abs}[1]{\lvert#1\rvert} 
\newcommand{\inner}[1]{\langle#1\rangle}

\theoremstyle{plain}
\newtheorem{theorem}{Theorem}[section]

\theoremstyle{definition}

\theoremstyle{remark}
\newtheorem{remark}[theorem]{Remark}

\newcommand{\MSC}[1]{{\narrower\noindent 2020 Mathematics Subject Classification}: #1\\}
\newcommand{\KEY}[1]{{\narrower\noindent Keywords: #1\\}}

\begin{document}

\title[Reconstruction for EIT using triangular decompositions]{Linearization-based direct reconstruction for EIT using triangular Zernike decompositions}

\author[A.~Autio]{Antti Autio}
\address[A.~Autio]{Department of Mathematics and Systems Analysis, Aalto University, P.O. Box~11100, 00076 Helsinki, Finland.}
\email{antti.2.autio@aalto.fi}

\author[H.~Garde]{Henrik Garde}
\address[H.~Garde]{Department of Mathematics, Aarhus University, Ny Munkegade 118, 8000 Aarhus C, Denmark.}
\email{garde@math.au.dk}

\author[M.~Hirvensalo]{Markus Hirvensalo}
\address[M.~Hirvensalo]{Department of Mathematics and Systems Analysis, Aalto University, P.O. Box~11100, 00076 Helsinki, Finland.}
\email{markus.hirvensalo@aalto.fi}

\author[N.~Hyv\"onen]{Nuutti Hyv\"onen}
\address[N.~Hyv\"onen]{Department of Mathematics and Systems Analysis, Aalto University, P.O. Box~11100, 00076 Helsinki, Finland.}
\email{nuutti.hyvonen@aalto.fi}

\begin{abstract}
This work implements and numerically tests the direct reconstruction algorithm introduced in [Garde \& Hyv\"onen, SIAM J.\ Math.\ Anal., 2024] for two-dimensional linearized electrical impedance tomography. Although the algorithm was originally designed for a linearized setting, we numerically demonstrate its functionality when the input data is the corresponding change in the current-to-voltage boundary operator. Both idealized continuum model and practical complete electrode model measurements are considered in the numerical studies, with the examined domain being either the unit disk or a convex polygon. Special attention is paid to regularizing the algorithm and its connections to the singular value decomposition of a truncated linearized forward map, as well as to the explicit triangular structures originating from the properties of the employed Zernike polynomial basis for the conductivity.
\end{abstract}

\maketitle

\KEY{electrical impedance tomography, linearization, direct reconstruction, Zernike polynomials, singular value decomposition}

\MSC{65N21, 65N20, 35R30, 35R25}

\section{Introduction} \label{sec:intro}

In {\em electrical impedance tomography} (EIT) the aim is to gain information on the conductivity distribution inside a physical body from current and potential measurements on its boundary; see the review articles \cite{Borcea02,Cheney99,Uhlmann09} and the references therein for more information on EIT. The reconstruction problem of EIT can be mathematically formulated as the ill-posed task of inverting the nonlinear operator that maps the internal conductivity of the examined body to the boundary measurements. 

The main goal of this paper is to implement and numerically test the direct reconstruction algorithm for two-dimensional EIT introduced in~\cite[Theorem~1.3]{Garde22}, capable of rapidly reconstructing any (even unbounded) $L^2$-perturbation to a homogeneous background conductivity. The article~\cite{Garde22} only considered the idealized {\em continuum model} (CM) of EIT that allows current input and potential measurement everywhere on the object boundary. In this work, we also test the algorithm on data simulated by resorting to the {\em complete electrode model} (CEM)~\cite{Cheng89,Somersalo92} as well as on real-world data from a water tank \cite{Hauptmann17}. Moreover, even though the reconstruction algorithm was originally designed for a linearized version of EIT, our numerical studies are exclusively based on difference data corresponding to the nonlinear forward map of either the CM or the CEM. The aim is to demonstrate that an algorithm based on a linearization can systematically produce good quality reconstructions in EIT.

The basic version of the algorithm is formulated assuming the examined object is modeled as the unit disk $D$, but it can also be applied to any smooth enough simply-connected domain by resorting to the Riemann mapping theorem, as explicitly outlined in the paragraph following~\cite[Theorem~1.3]{Garde22}. In the framework of the linearized CM, one assumes that the available data is the Fr\'echet derivative of the forward map, evaluated at the unit conductivity, in the direction of the to-be-reconstructed conductivity perturbation. This is a linear operator sending input currents to measured boundary potentials. To be more precise, the algorithm takes as its input the elements in the infinite matrix representation of this boundary map with respect to the standard Fourier basis. As hinted above, in our numerical studies these matrix elements are replaced by those corresponding to the actual change in the {\em Neumann-to-Dirichlet} (ND) map when the conductivity inside $D$ is perturbed. The reconstruction algorithm itself is comprised of a set of explicit linear formulas that map the aforementioned matrix elements to the coefficients, $c_{j,k}$, in the expansion of the conductivity perturbation in the Zernike polynomial basis~\cite{Zernike1934}; see also \cite{Allers1991} for a similar approach.

The Zernike polynomials are indexed by two indices: $j \in \Z$ controlling the Fourier frequency in the polar angle and $k \in \N_0$ controlling (for a fixed $j$) the behavior in the radial variable, with increasing $k$ indicating more oscillations closer to the origin. As explicitly proven in~\cite{Garde22}, reconstructing the coefficients $c_{j,k}$ becomes more unstable as $k$ increases. On the other hand, it turns out that $c_{j,k}$ with a fixed angular index $j$ only affect the $j$th diagonal in the matrix representation of the associated Fr\'echet derivative. In consequence, a certain angular frequency in the conductivity perturbation can be reconstructed by only considering a single diagonal in the data matrix, and the whole reconstruction algorithm can thus be decoupled into independent smaller tasks indexed by $j$. What is more, the linear dependence of the $j$th diagonal in the data matrix on $c_{j,k}$ can be represented by an infinite triangular matrix. The formulas comprising the direct reconstruction algorithm correspond to explicitly inverting these triangular matrices for all $j \in \Z$. Note that singular vectors of the linearized forward map must also obey such a decoupling, that is, a right-hand singular vector may only contain a single angular frequency and a left-hand singular vector can only have nonzero elements on a single diagonal of its matrix representation.

The decoupling of the angular frequencies allows, in particular, to separately regularize the subsystems corresponding to different angular indices $j$. This leads to increasing level of regularization as a function of $\abs{j}$. Numerically considering radial indices beyond, say, $k=15$, is impractical due to the increasing ill-posedness of the subproblems for higher radial indices. Hence, any reasonable regularization method can be used for the decoupled subproblems without loosing the instantaneous nature of the reconstruction algorithm. We consider two simple approaches accompanied by the Morozov discrepancy principle: 
\begin{enumerate}[(i)]
    \item truncated {\em singular value decomposition} (SVD) and
    \item choosing the sizes of the decoupled subsystems as a function of $\abs{j}$ by monitoring the absolute values of the elements on their diagonals.
\end{enumerate}
The former has a better theoretical justification, but it cannot fully exploit the triangular structure of the decoupled systems as can the latter.

\subsection{Article structure}

Section \ref{sec:linearization} introduces the linearized CM for the unit disk and extends it for square-integrable conductivity perturbations following~\cite{Garde22}. In Section \ref{sec:algorithm}, the linearized forward map is written out explicitly with respect to the Zernike polynomial basis for the conductivity perturbation and the Fourier basis for the boundary measurements. In addition, the triangular structures of the subsystems corresponding to different angular Zernike indices are studied. The implementation of regularized reconstruction algorithms are considered in Section~\ref{sec:implementation}, while Section~\ref{sec:poly_CEM} explains how data for the algorithm can be measured/simulated on polygonal domains or based on the CEM. The numerical experiments are presented in Section~\ref{sec:res}.

\section{Linearization of the CM of EIT} \label{sec:linearization}

This section introduces the linearization of the CM with respect to square-integrable conductivity perturbations following the material in~\cite{Garde22}. In particular, our approach only requires considering linearizations for the CM in the unit disk $D \subset \R^2$. In the numerical experiments with polygonal domains or the CEM, we rely on methods for replicating or approximating CM boundary data for the unit disk based on (electrode) measurements on other simply-connected planar domains; see \cite{Garde21,Hyvonen17d}.

Consider an isotropic conductivity
\begin{equation*}
    \gamma \in L^\infty_+(D) \coloneqq \big \{ \kappa \in L^\infty(D) :  \text{ess} \inf(\redel \kappa) > 0 \big \}
\end{equation*}
in $D$ and a normal current density 
\begin{equation*}
    f \in L^2_\diamond(\partial D) \coloneqq \big \{ g \in L^2(\partial D) : \inner{g, 1}_{L^2(\partial D)} = 0 \big \}
\end{equation*}
on $\partial D$. According to the CM, the corresponding electric potential $u$ in $D$ weakly satisfies a Neumann problem for the conductivity equation:
\begin{equation} \label{eq:strong_form}
    - \nabla \cdot (\gamma \nabla u) = 0 \quad \text{in } D, \qquad  \nu \cdot (\gamma \nabla u) = f \quad \text{on } \partial D,
\end{equation}
where $\nu$ is the exterior unit normal of $\partial D$.
The variational form of \eqref{eq:strong_form} is to find $u \in H^1(D)$ that satisfies
\begin{equation} \label{eq:weak_form}
    \int_D \gamma \nabla u \cdot \overline{\nabla v} \, \di x =  \inner{f, v}_{L^2(\partial D)} \qquad \text{for all } v \in H^1(D),
\end{equation}
for which the Lax--Milgram lemma yields a unique solution
\begin{equation*}
    u^\gamma_f \in H^1_\diamond(D) \coloneqq \{ w \in H^1(D) \, : \, \inner{w|_{\partial D}, 1}_{L^2(\partial D)} = 0 \}.
\end{equation*}
Note that the use of the subspace $H^1_\diamond(D) \subset H^1(D)$ corresponds to a systematic way of choosing the ground level of potential.

The associated ND, or current-to-voltage, map 
\begin{equation} \label{eq:ND_map}
    \Lambda(\gamma):  f \mapsto  u^\gamma_f|_{\partial D}
\end{equation}
is a compact self-adjoint operator in \( \mathscr{L}(L^2_\diamond(\partial D)) \). In fact, in the considered two-dimensional setting with a smooth domain,  $\Lambda(\gamma)$ belongs to the space \( \mathscr{L}_{\rm HS}(L^2_\diamond(\partial D)) \) of Hilbert--Schmidt operators. The forward map of the CM is defined as
\begin{equation*}
    \Lambda:
    \left\{
    \begin{array}{l}
      \gamma \mapsto \Lambda(\gamma), \\[1mm]
       L^\infty_+(D) \to \mathscr{L}_{\rm HS}(L^2_\diamond(\partial D)) \subset \mathscr{L}(L^2_\diamond(\partial D)).
    \end{array}
    \right.
\end{equation*}
It is well-known that $\Lambda$ is Fr\'echet differentiable at any $\gamma  \in L^\infty_+(D)$, with the corresponding derivative $D \Lambda (\gamma)$ in $\mathscr{L}( L^\infty(D), \mathscr{L}_{\rm HS}(L^2_\diamond(\partial D) ) ) \subset \mathscr{L}( L^\infty(D), \mathscr{L}(L^2_\diamond(\partial D) ) )$; see,~e.g.,~\cite{garde2021series}.

In this work, we are interested in the Fr\'echet derivative of $\Lambda$ at the unit conductivity,~i.e.,~$F = D \Lambda(1)$. $F$ satisfies
\begin{equation} \label{eq:frechet_redef}
    \inner{(F \eta) f, g}_{L^2(\partial D)} = - \int_D \eta \nabla u_f^1 \cdot \overline{\nabla u_g^1} \, \di x
\end{equation}
for any $\eta \in L^\infty(D)$ and all $f, g \in L^2_\diamond(\partial D)$. In fact, \eqref{eq:frechet_redef} uniquely characterizes $F$ as an element of $\mathscr{L}( L^\infty(D), \mathscr{L}(L^2_\diamond(\partial D) ) )$, and it can thus be used as a definition.

According to the elliptic regularity theory, the gradients of the harmonic functions $u_f^1$ and $u_g^1$ belong to $H^{1/2}(D)$, with their norms, respectively, bounded by the $L^2(\partial D)$-norms of $f$ and $g$. By the continuous embedding $H^{1/2}(D) \hookrightarrow L^4(D)$, the formula \eqref{eq:frechet_redef} can thus be used to extend $F$ to an operator in $\mathscr{L} ( L^2(D), \mathscr{L}(L^2_\diamond(\partial D) ) )$. This is how we interpret $F$ in what follows; consult \cite{Garde22} for the details of this extension argument.

\section{Explicit representation for $F$ and its inverse} \label{sec:algorithm}

In this section, we continue to assume that the investigated domain is the unit disk $D$ and give an explicit representation for the bounded linear map
\begin{equation*}
    F: \left\{
    \begin{array}{l}
        \eta \mapsto  D \Lambda(1) \eta, \\[1mm]
        L^2(D) \to \mathscr{L}(L^2_\diamond(\partial D))
    \end{array}
    \right.
\end{equation*}
with respect to the Zernike polynomial basis~\cite{Zernike1934} for $L^2(D)$ and the standard Fourier basis for $L^2_\diamond(\partial D)$. In fact, we also present an explicit characterization for the inverse of $F$ on its range. Take note that the injectivity of $F$ follows from the original uniqueness proof of Calder\'on~\cite{Calderon80} that generalizes for square-integrable conductivity perturbations in two dimensions; see also \cite[Section~6]{Garde22} for an alternative proof. We still partially follow the material and presentation in~\cite{Garde22}.
    
\subsection{Zernike polynomials and Fourier basis} \label{subsec:Zernike}

As in \cite{Garde22}, we define the Zernike polynomials in the polar coordinates $(r,\theta)$ as 
\begin{equation} \label{eq:zernike}
    \psi_{j, k}(r, \theta) = \sqrt{\frac{\abs{j} + 2 k + 1}{\pi}} R^{\abs{j}}_{\abs{j} + 2 k}(r) \, \e^{\I j\theta}, \qquad j \in \Z, \ k \in \N_0,
\end{equation}
where 
\begin{equation*}
    R^{\abs{j}}_{\abs{j} + 2 k}(r) = \sum^k_{i=0} (-1)^i \binom{\abs{j} + 2 k - i}{i} \binom{\abs{j} + 2 k - 2 i}{k - i} \, r^{\abs{j} + 2 k - 2 i}
\end{equation*}
is a Zernike radial function. As mentioned in Section~\ref{sec:intro}, the two indices  \( j \in \Z \) and \( k \in \N_0 \) control the angular and radial characteristics of the Zernike polynomials, with larger absolute values of $j$ and $k$ indicating faster spatial oscillations.

When scaled as in \eqref{eq:zernike}, the Zernike polynomials form an orthonormal basis for $L^2(D)$, and they are used for representing the conductivity perturbation $\eta$ as
\begin{equation} \label{eq:Zexpansion}
    \eta = \sum_{j=-\infty}^\infty \eta_j = \sum_{j= -\infty}^\infty \sum_{k=0}^\infty c_{j,k} \psi_{j,k},
\end{equation}
where we call $\eta_j$ the \emph{$j$th angular frequency in $\eta$}, defined by
\begin{equation*}
    \eta_j = \sum_{k=0}^\infty c_{j,k} \psi_{j,k}, \qquad j \in \Z.
\end{equation*}
In the following, we implicitly assume such a Zernike polynomial expansion with square-summable coefficients $\{c_{j,k}\}_{j\in \Z, k \in \N_0} \subset \C$ for the to-be-reconstructed conductivity perturbation~$\eta$.

The employed orthonormal basis for $L^2_\diamond(\partial D)$ is a standard Fourier basis
\begin{equation*}
    f_m(\theta) = \frac{1}{\sqrt{2 \pi}}\e^{\I m\theta}, \qquad m \in \Z',
\end{equation*}
where $\Z' = \Z \setminus \{ 0 \}$ in accordance with all elements of $L^2_\diamond(\partial D)$ having zero mean. With these choices, the mapping $F: L^2(D) \to \mathscr{L}(L^2_\diamond(\partial D))$ is completely characterized by the coefficients
\begin{equation} \label{eq:a_coef}
    a^{j,k}_{m,n} = \inner{(F \psi_{j,k}) f_m, f_n}_{L^2(\partial D)}, \qquad j \in \Z, \ k \in \N_0, \ m,n \in \Z'.
\end{equation}
Indeed, for a given conductivity perturbation $\eta$ expanded as in \eqref{eq:Zexpansion}, the mapping $F\eta \in \mathscr{L}(L^2_\diamond(\partial \Omega))$ is represented in the Fourier basis as an infinite matrix given element-wise as
\begin{equation} \label{eq:inf_matrix}
     a_{m,n}(\eta) = \inner{(F \eta) f_m, f_n}_{L^2(\partial D)} = \sum_{j= -\infty}^\infty a_{m,n}(\eta_j) = \sum_{j= -\infty}^\infty \sum_{k=0}^\infty a^{j,k}_{m,n} c_{j,k},   \qquad m,n \in \Z',
\end{equation}
where
\begin{equation} \label{eq:j_matrix}
   a_{m,n}(\eta_j) = \sum_{k=0}^\infty a^{j,k}_{m,n} c_{j,k}, \qquad j \in \Z, \  m,n \in \Z', 
\end{equation}
only depends on the Zernike coefficients defining $\eta_j$. In particular, $\{a_{m,n}(\eta) \}_{m,n \in \Z'}$ can be considered the data for the linearized inverse problem of the~CM in the unit disk.

\subsection{Matrix representation of the linearized model}

Our reconstruction algorithm is based on the fundamental observation that most matrix elements in \eqref{eq:a_coef} vanish. Indeed, \cite[Eq.~(4.5)]{Garde22} immediately gives
\begin{equation} \label{eq:expl_coef}
    a^{j,k}_{m,n} = - \frac{1}{\sqrt{\pi}} \, \frac{\sqrt{\abs{j} + 2 k + 1}}{\min \{\abs{m},\abs{n} \} + \abs{j} + k} \, \prod^k_{i=1} \frac{\min \{\abs{m},\abs{n} \} - i}{\abs{j} + \min \{\abs{m},\abs{n} \} + k - i}
\end{equation}
if $n=m+j$, $k < \min \{\abs{m},\abs{n} \}$ and $mn > 0$, and $a^{j,k}_{m,n}= 0$ for all other combinations of $j \in \Z$, $k \in \N_0$ and $m,n \in \Z'$. When $k=0$, the product in \eqref{eq:expl_coef} is defined to take the value $1$. In particular, if $\{ a^{j,k}_{m,n} \}_{n,m \in \Z'}$ are interpreted as elements of an infinite matrix (cf.~\eqref{eq:inf_matrix}), then for a fixed $j$ but any $k$, all nonzero elements in the matrix lie on its $j$th diagonal. This means that the infinite matrix defined by \eqref{eq:j_matrix} also only has nonzero elements on its $j$th diagonal, so by linearity we have
\begin{equation} \label{eq:jdiagreduction}
    a_{m,m+j}(\eta) = a_{m,m+j}(\eta_j), \qquad j\in \Z,\ m\in\Z'. 
\end{equation}
In other words, the $j$th angular frequency $\eta_j$ in the conductivity perturbation,~i.e.,~the Zernike coefficients $\{ c_{j,k} \}_{k\in \N_0}$, only contributes to the $j$th diagonal in the matrix representation \eqref{eq:inf_matrix} of $F\eta$, which is the decoupling in the angular frequencies we were looking for.

To explain how this  decoupling can be exploited in solving for the coefficients $\{ c_{j,k} \}_{j \in \Z, k \in \N_0}$ of $\eta$ from the data $\{ a_{m,n}(\eta) \}_{m,n \in \Z'}$, let us define two infinite vectors indexed by the angular index $j \in \Z$:
\begin{equation} \label{eq:inf_vec}
    \mathbf{c}^j \coloneqq [c_{j,k-1}]_{k=1}^\infty, \qquad \mathbf{a}^{j} \coloneqq \left\{
    \begin{array}{ll}
        \big[ a_{m,m+j}(\eta) ]_{m=1}^{\infty} & \quad \text{if } j \geq 0, \\[2mm]
        \big[ a_{-m,-m+j}(\eta) ]_{m=1}^{\infty} & \quad \text{if } j < 0.
    \end{array}
    \right.
\end{equation}
The Zernike coefficients defining $\eta_j$ form $\mathbf{c}^j$, whereas $\mathbf{a}^{j}$ is composed of certain elements on the $j$th diagonal of $\{ a_{m,n}(\eta) \}_{m,n \in \Z'}$. Our aim is to present the linear relationship between $\mathbf{a}^{j}$ and $\mathbf{c}^j$ as a  multiplication by an infinite lower triangular matrix.

Let us start with the case $j \geq 0$. According to \eqref{eq:jdiagreduction} and \eqref{eq:j_matrix}, and since $a^{j,k}_{m,n} = 0$ unless $n=m+j$, $k < \min \{\abs{m},\abs{n} \}$ and $mn > 0$,
\begin{equation} \label{eq:aplus}
    \mathbf{a}^{j}_m = \sum_{k=0}^\infty a^{j,k}_{m,m+j} c_{j,k} = \sum_{k=0}^{m-1} a^{j,k}_{m,m+j} c_{j,k} = \sum_{k=1}^{m} a^{j,k-1}_{m,m+j} \mathbf{c}^j_{k}, \qquad m \in \N,
\end{equation}
which is the sought-for infinite triangular system for $j \geq 0$. Analogously, for $j < 0$,
\begin{equation} \label{eq:aminus}
    \mathbf{a}^{j}_m =  \sum_{k=0}^\infty a^{j,k}_{-m,-m+j} c_{j,k} = \sum_{k=0}^{m-1} a^{j,k}_{-m,-m+j} c_{j,k} = \sum_{k=1}^{m} a^{j,k-1}_{-m,-m+j} \mathbf{c}^j_{k}, \qquad m \in \N.
\end{equation}
Due to obvious symmetries in \eqref{eq:expl_coef}, $a^{-\abs{j},k-1}_{-m,-m-\abs{j}} = a^{\abs{j},k-1}_{m,m+\abs{j}}$, and thus the coefficients of the triangular systems in \eqref{eq:aplus} and \eqref{eq:aminus} only depend on $j$ via $\abs{j}$. Hence,
\begin{align} \label{eq:triang_syst}
    F^{\abs{j}} \mathbf{c}^j &= \mathbf{a}^{j}, \qquad j \in \Z,
\end{align}
where the $\abs{j}$-dependent infinite lower triangular  system matrix is defined by
\begin{equation} \label{eq:inf_F_matrix}
    F^{\abs{j}}_{m,k} = \left\{
    \begin{array}{ll}
        a^{\abs{j},k-1}_{m,m+\abs{j}} \quad & \text{if } 1 \leq k  \leq m \in \N, \\[1mm]
        0 \quad & \text{otherwise}.
    \end{array}
    \right.
\end{equation}
In particular, note that the diagonal elements $F^{\abs{j}}_{m,m}$ are nonzero for all $m \in \N$ and $j \in \Z$, which means that one can, in principle, solve for the components of $\mathbf{c}^j$  via forward substitution.

\begin{remark}
    Observe that less than a half of the $j$th diagonal of $\{ a_{m,n}(\eta) \}_{n,m\in \Z'}$ is included in  $\mathbf{a}^{j}$. As this is already sufficient for uniquely determining $\eta_j$, the data that is left out must be redundant. Indeed, it follows straightforwardly from \eqref{eq:inf_matrix} and \eqref{eq:expl_coef} that the $j$th diagonal of $\{ a_{m,n}(\eta) \}_{n,m\in \Z'}$ is symmetric with respect to its midpoint,~i.e.,~the (possibly virtual) element $a_{-j/2,j/2}(\eta)$. Moreover, the $\abs{j}$ central elements on the $j$th diagonal vanish since their indices have opposite signs (see the text after \eqref{eq:expl_coef}). In fact, two quadrants of $\{ a_{m,n}(\eta) \}_{n,m\in \Z'}$ are zeros, and either of the two remaining ones carry all information contained in the matrix.

    In particular, one could equivalently define the $j$th data vector by employing the opposite ends of the diagonals:
    \begin{equation} \label{eq:new_inf_vec}
        \mathbf{a}^{j} = \left\{
        \begin{array}{ll}
            \big[ a_{-m-j,-m}(\eta) ]_{m=1}^{\infty} & \quad \text{if } j \geq 0, \\[2mm]
            \big[ a_{m-j,m}(\eta) ]_{m=1}^{\infty} & \quad \text{if } j < 0.
        \end{array}
        \right.
    \end{equation}
    If the whole matrix $\{ a_{m,n}(\eta) \}_{n,m\in \Z'}$ is available, this observation makes it possible to slightly reduce the effect of measurement noise by using the average of the two equivalent versions of $\mathbf{a}^j$ in \eqref{eq:inf_vec} and \eqref{eq:new_inf_vec} as the data in numerical considerations.
\end{remark}

\subsection{Inversion of the infinite triangular systems}

The angularly decoupled triangular systems in \eqref{eq:triang_syst} imply a direct algorithm for reconstructing the coefficients $c_{j,k}$ one angular frequency at a time, by resorting to simple forward substitution. Since the resulting formulas have already been worked out in \cite[Theorem~1.3]{Garde22} --- although without emphasizing the underlying triangular structure or identifying the method as forward substitution ---, we settle for writing the resulting reconstruction formulas in our notation:
\begin{align} \label{eq:recursion}
    \mathbf{c}^j_{k+1} &= -\sqrt{\pi (\abs{j} + 2k + 1)}  \; \binom{\abs{j} + 2k}{k} \; \mathbf{a}^j_{k+1} \notag\\
    &\phantom{{}=}-\sum^{k}_{i=1} \mathbf{c}^j_i \frac{\sqrt{(\abs{j} + 2k + 1) (\abs{j} + 2i - 1)}}{\abs{j} + k + i} \, \binom{\abs{j} + 2k}{k - i + 1}, \qquad k= 0,1,2,\dots \, , 
\end{align}
for any $j \in \Z$. When $k=0$, the sum on the right-hand side of \eqref{eq:recursion} is interpreted to vanish.

By virtue of the definition of $\mathbf{a}^j$ via \eqref{eq:inf_matrix} and \eqref{eq:inf_vec}, one only needs to apply two boundary currents $f_{\pm 1}$ and measure (all Fourier frequencies in) the corresponding potentials to reconstruct $c_{j,0} = \mathbf{c}^{j}_1$ for all $j \in \Z$. Moving to a next $k$-level always `costs' applying two more Fourier currents in increasing order with respect to the absolute value of their angular frequencies. On the other hand, if one resorts to the dual definition of $\mathbf{a}^j$ in \eqref{eq:new_inf_vec}, reconstructing $c_{j,0}$ for all $j \in \Z$ requires applying all Fourier boundary currents but only recording the two lowest Fourier modes in the resulting boundary potentials. In this case, moving to a next $k$-level always requires recording two more Fourier frequencies in all measured boundary potentials.

\section{Regularization for truncated systems} \label{sec:implementation}

This section discusses the implementation of regularized reconstruction algorithms based on truncations of the triangular systems~\eqref{eq:triang_syst}. To this end, we introduce the notation 
\begin{equation*}
    \Z_M \coloneqq \{-M,\dots,M\} \quad \textup{and} \quad \Z'_M \coloneqq \Z_M\setminus\{0\}
\end{equation*}
and assume that the available data is the truncated matrix $\{ a_{m,n} \}_{m,n \in \Z'_M}$ for some fixed $M\in\N$. This means that we assume to be able to apply the Fourier boundary currents indexed by $\Z'_M$ and record the corresponding Fourier frequencies in the resulting relative boundary potentials. As a consequence, the infinite linear systems~\eqref{eq:triang_syst} are to be replaced by finite-dimensional ones due to a shortage of data:
\begin{align} \label{eq:finite_triang_syst}
     F^{\abs{j},M} \mathbf{c}^{j,M} &= \mathbf{a}^{j,M},
\end{align}
where $F^{\abs{j},M}$ is a lower triangular square matrix and $\mathbf{a}^{j,M}$ and $\mathbf{c}^{j,M}$ are, respectively, composed of a certain number of first components in $\mathbf{a}^{j}$ and $\mathbf{c}^{j}$ defined by~\eqref{eq:inf_vec}.

The length of the vectors $\mathbf{a}^{j,M}$ vary as a function of $j$ as they correspond to different diagonals of a {\em finite} square matrix. Indeed, $\mathbf{a}^{j,M} \in \C^{M-\abs{j}}$, which also reveals that the systems \eqref{eq:finite_triang_syst} are nonempty only for $j \in \Z_{M-1}$. To summarize,
\begin{equation*}
    F^{\abs{j},M} \in \R^{(M-\abs{j}) \times (M-\abs{j})}, \quad \mathbf{a}^{j,M}, \mathbf{c}^{j,M} \in \C^{M-\abs{j}}, \qquad j \in \Z_{M-1},
\end{equation*}
defined, respectively, by the top left corner of $F^{\abs{j}}$ in \eqref{eq:inf_F_matrix} and the first $M-\abs{j}$ components of the vectors in~\eqref{eq:inf_vec}. Thanks to the triangular structure of \eqref{eq:finite_triang_syst}, the truncated systems~\eqref{eq:finite_triang_syst} are still uniquely solvable. In the framework of the linearized CM, they are also exact in the sense that the components of the solution vectors $\mathbf{c}^{j,M}$ are precisely the corresponding Zernike coefficients of the to-be-reconstructed conductivity perturbation $\eta$. In particular, the explicit inversion formulas \eqref{eq:recursion} are still valid, but only for $j \in \Z_{M-1}$ and $k = 0, \dots, M-\abs{j}-1$.

\subsection{Truncated SVD} \label{sec:svd_regularization}

\begin{figure}[t]
    \centering
    \includegraphics[width=.5\textwidth]{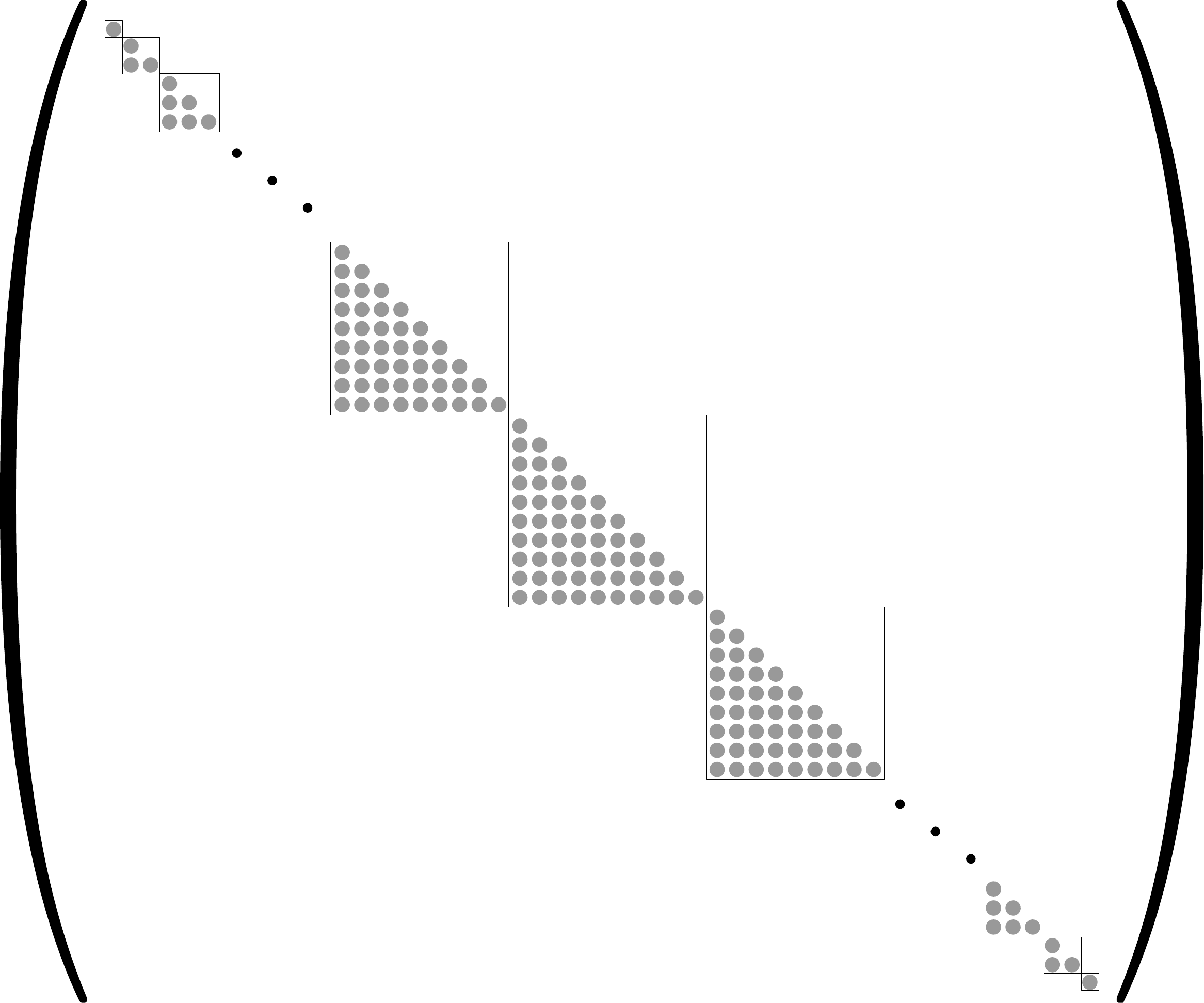}
    \caption{The structure of $F^M$, for $M=10$, as a block diagonal matrix with lower triangular blocks.}
    \label{fig:matstructure}
\end{figure}

The first considered regularization technique for the $2 M - 1$ triangular systems \eqref{eq:finite_triang_syst} is truncated SVD combined with the Morozov discrepancy principle. Let us denote by $\mathbf{a}^M, \mathbf{c}^M \in \C^{M^2}$, respectively, the concatenations of the data vectors $\mathbf{a}^{j,M} \in \C^{M-\abs{j}}$ and Zernike coefficient vectors $\mathbf{c}^{j,M} \in \C^{M-\abs{j}}$ for $j \in \Z_{M-1}$. The systems \eqref{eq:finite_triang_syst} can be written as a single matrix equation
\begin{equation} \label{eq:finite_complete_system}
    F^M \mathbf{c}^{M} = \mathbf{a}^M,
\end{equation}
where
\begin{equation} \label{eq:bloc_triangular}
    F^M = {\rm diag}\big(F^{M-1,M}, F^{M-2,M}, \dots, F^{0,M}, \dots, F^{M-2,M}, F^{M-1,M} \big) \in \R^{M^2 \times M^2}
\end{equation}
is a block diagonal matrix whose structure is visualized in Fig.~\ref{fig:matstructure} for $M=10$.
 
As a consequence, singular systems for the individual $F^{l,M}$, $l=0,\dots, M-1$, completely define an SVD for the whole $F^M$. Indeed, each singular value of any $F^{l,M}$, $l=0, \dots, M-1$, is also a singular value for $F^M$, with the corresponding singular vectors of $F^M$ being suitably zero-padded versions of those of $F^{l,M}$. Moreover, excluding the case $l=0$, each singular value of $F^{l,M}$, $l=1,\dots, M-1$, appears twice in an SVD for $F^M$ since $F^{l,M}$ also appears twice on the diagonal of $F^M$.

Computing an SVD for $F^M$ thus boils down to computing such for the $M$ small matrices $F^{l,M}$, $l=0, \dots, M-1$,~i.e.,~for slightly more than a half of the diagonal blocks defining $F^M$. This is computationally tractable in practice. However, to enable adopting truncated SVD as a regularization method, the singular values of the complete system matrix $F^M$ also need to be ordered. For the sake of symmetry, we always include in a truncated SVD of $F^M$ either both or neither of the two right-left singular vector pairs of $F^M$ associated to the same singular value of a repeated block $F^{l,M}$, $l=1, \dots, M-1$, in \eqref{eq:bloc_triangular}. This leads to $M(M+1)/2$ potentially different singular values for $F^M$. The ordering of these $M(M+1)/2$ singular values is then achieved with the help of a mapping
\begin{equation*}
    \tau: \{1, 2, \dots, M(M+1)/2 \} \to \N_0^{M}
\end{equation*}
that works as follows: $\tau_l(p)$,~i.e.,~the $l$th component of $\tau(p)$, indicates how many of the largest $p$ singular values of $F^M$ correspond to $F^{l,M}$. The map $\tau$ can be formed and tabulated offline if $M$ is known prior to taking the measurements.
 
Let $(F^{\abs{j},M})^\dagger_q$ denote a truncated Moore--Penrose pseudoinverse of $F^{\abs{j},M}$,~i.e., the pseudoinverse that has been formed after all but the largest $q \in \{0, \dots, M-\abs{j}\}$ singular values of $F^{\abs{j},M}$ have been replaced by zeros in an SVD for $F^{\abs{j},M}$. For $p \in \{1, 2, \dots, M(M+1)/2 \}$, acting as a truncation index/regularization parameter, we define the truncated SVD solution for \eqref{eq:finite_triang_syst} via
\begin{equation*}
    \mathbf{c}^{j,M,p} = (F^{\abs{j},M})^\dagger_{\tau_{\abs{j}}(p)} \mathbf{a}^{j,M} , \qquad j \in \Z_{M-1}.
\end{equation*}
The complete truncated SVD solution $\mathbf{c}^{M,p} \in \C^{M^2}$ for \eqref{eq:finite_complete_system} is then defined as the concatenation of $\mathbf{c}^{j,M,p}$ for $j \in \Z_{M-1}$.

Assume now that the available data is a noisy version $\mathbf{a}^M_\delta$ of $\mathbf{a}^M$, approximately satisfying
\begin{equation} \label{eq:noisy_data}
    \norm{\mathbf{a}^M_\delta - \mathbf{a}^M}_2 = \delta > 0,
\end{equation}
where $\norm{\, \cdot \,}_2$ denotes the Euclidean norm.\footnote{Employing the Euclidean norm for the diagonals of the data matrix corresponds to the Frobenius norm for the matrix itself, which is in line with the Hilbert--Schmidt topology for the ND map in~\cite{Garde22}.} Let $\mathbf{c}^{M,p}_{\delta}$ be a truncated SVD solution corresponding to the noisy data $\mathbf{a}^M_\delta$. Truncated SVD solutions are combined with the Morozov discrepancy principle in the standard manner,~i.e.,~by choosing the smallest truncation index $p \in \{1, 2, \dots, M(M+1)/2 \}$ that satisfies
\begin{equation}
\label{eq:morozov1}
    \norm{F^M \mathbf{c}^{M,p}_{\delta} - \mathbf{a}^M_\delta}_2 \leq \omega \delta,
\end{equation}
where $\omega \geq 1$ is a `fudge' factor that aims to account for linearization and numerical errors. The choice of $\omega$ is considered in Section~\ref{sec:res}.

\subsection{Truncation of the triangular subsystems} \label{sec:triangular_regularization}

A downside of truncated SVD is that it does  not fully exploit the triangular structure of~\eqref{eq:finite_triang_syst} and, in particular, the explicit solution formula~\eqref{eq:recursion}. Our second regularization approach aims to address this shortcoming by directly resorting to~\eqref{eq:recursion}, but only up to a $\abs{j}$-dependent threshold for the radial index $k$.

It is straightforward to prove that the diagonal elements in each $F^{j, M}$, $j=0, \dots, M-1$, are monotonically decreasing in their absolute values. Running the recursion~\eqref{eq:recursion} until, say, $k+1 = q \leq M - \abs{j}$ corresponds to applying forward substitutions to operate on the first $q$ elements in $\mathbf{a}^{j, M}$ with the inverse of the triangular square matrix composed of the first $q$ rows and columns in $F^{\abs{j}, M}$. The most unstable single step in this procedure is thus dividing by the last involved diagonal element $F^{\abs{j}, M}_{q,q}$. Motivated by this observation, our second (more heuristic) regularization strategy is based on satisfying the Morozov discrepancy principle by running~\eqref{eq:recursion} for all $j \in \Z_{M-1}$ in such a way that the smallest absolute value of a diagonal element involved in a forward substitution is maximized. Although such an idea does not fully account for the accumulation of errors in forward substitutions, the smallest absolute value of a diagonal element in a lower triangular matrix is in any case closely related to the stability of forward substitutions for solving a corresponding linear system; see~\cite[Chapter~8]{Higham02} for more information.

To simplify our notation, let $(F^{\abs{j},M})^{-1}_q$ denote the inverse of the top left triangular $q \times q$ block in $F^{\abs{j},M}$ and let $\mathbf{a}^{j, M,q}$ be the vector composed of the first $q$ elements in $\mathbf{a}^{j, M}$. Applying $(F^{\abs{j},M})^{-1}_q$ to $\mathbf{a}^{j, M,q}$ is equivalent to running the algorithm~\eqref{eq:recursion} until $k+1 = q$ for the considered $j \in  \Z_{M-1}$, with $q \in \{0, \dots, M - \abs{j}\}$. Note that $F^M$ contains $M(M+1)/2$ potentially distinct diagonal elements, and let
\begin{equation*}
    \upsilon: \{1, 2, \dots, M(M+1)/2 \} \to \N_0^{M}
\end{equation*}
order the absolute values of those elements in the same way that $\tau$ orders singular values: $\upsilon_l(p)$ indicates how many of the in-absolute-value largest $p$ diagonal elements of $F^M$  are found on the diagonal of $F^{l,M}$. Note that $\upsilon$ can be formed and tabulated offline without any prior knowledge apart from some upper bound on $M$.

Given $p \in \{1, 2, \dots, M(M+1)/2 \}$, the {\em truncated triangular solution} to \eqref{eq:finite_triang_syst} is defined via
\begin{equation} \label{eq:tts}
    \mathbf{c}^{j,M,p} = \begin{bmatrix} (F^{\abs{j},M})^{-1}_{\upsilon_{\abs{j}}(p)} \mathbf{a}^{j,M,\upsilon_{\abs{j}}(p)} \\[2mm] \mathbf{0} \end{bmatrix} , \qquad j \in \Z_{M-1},
\end{equation}
where $\mathbf{0}$ is a zero-vector of length $M - \abs{j} - \upsilon_{\abs{j}}(p)$. The complete truncated triangular solution $\mathbf{c}^{M,p} \in \C^{M^2}$ for \eqref{eq:finite_complete_system} is then defined as the concatenation of $\mathbf{c}^{j,M,p}$ for $j  \in \Z_{M-1}$.

Assume then that the available data is a noisy version $\mathbf{a}^{M}_\delta$ of $\mathbf{a}^M$, approximately satisfying \eqref{eq:noisy_data}, and let $\mathbf{c}^{M,p}_{\delta}$ be a corresponding truncated triangular solution. Truncated triangular solutions are combined with the Morozov discrepancy principle by choosing the smallest truncation index $p \in \{1, 2, \dots, M(M+1)/2 \}$ that satisfies
\begin{equation} \label{eq:morozov2}
    \norm{F^M \mathbf{c}^{M,p}_{\delta} - \mathbf{a}^M_\delta}_2^2 = \sum_{j \in \Z_{M-1}} \norm{\widetilde{\mathbf{a}}^{M,j, \upsilon_{\abs{j}}(p)}_{\delta}}_2^2 \leq \omega^2 \delta^2,
\end{equation}
where $\omega \geq 1$ is again a `fudge' factor and $\widetilde{\mathbf{a}}^{M,j, \upsilon_{\abs{j}}(p)}_{\delta}$ denotes the vector composed of the last $M-\abs{j}-\upsilon_{\abs{j}}(p)$ elements in $\mathbf{a}^{M,j}_\delta$. The equality in~\eqref{eq:morozov2} follows from the use of the exact inverse matrix for the top left $\upsilon_{\abs{j}}(p) \times \upsilon_{\abs{j}}(p)$ subsystem in~\eqref{eq:tts}.

\section{Simulating data on piecewise smooth domains and with the CEM} \label{sec:poly_CEM}

In this section, we explain how data measured on more general planar domains or via electrode measurements can be used as input for the regularized algorithms introduced in the preceding section. Although the algorithms have been designed for linearized measurements, all data used as their input in this work are obtained from boundary potential measurements corresponding to nonlinear forward models of EIT.

\subsection{Other planar domains with CM} \label{sec:piecewise_smooth}

Assume the imaged object is modeled as a simply-connected domain $\Omega \subset \R^2$ with a Lipschitz boundary that is composed of a finite number of $C^{1,\alpha}$ smooth arcs. Let $\Psi : \Omega \to D$ be a bijective conformal mapping whose existence is guaranteed by the Riemann mapping theorem. Denote the inverse of $\Psi$ by $\Phi$ and their complex derivatives by $\Psi'$ and $\Phi'$. The mappings $\Psi$ and $\Phi$ have continuously differentiable extensions to the relative interiors of the $C^{1,\alpha}$ arcs composing $\partial \Omega$ and their images under $\Psi$, respectively~\cite{Pommerenke92}. The absolute values $\abs{\Psi'}$ and $\abs{\Phi'}$ restricted onto $\partial \Omega$ and $\partial D$, respectively, are the Jacobian determinants of the corresponding boundary transformations.

On $\partial \Omega$ define
\begin{equation} \label{eq:conformalfm}
    \widetilde{f}_m \coloneqq (f_m \circ \Psi) \abs{\Psi'}, \qquad m \in \Z',
\end{equation}
with $f_m$ still denoting the $m$th Fourier basis function on $\partial D$. By performing a change of variables on each of the $C^{1, \alpha}$ arcs composing $\partial \Omega$, it follows that
\begin{equation*}
    \int_{\partial \Omega} \abs{\widetilde{f}_m}^2 \, \di s = \int_{\partial \Omega} \abs{f_m  \circ \Psi}^2 \abs{\Psi'}^2 \, \di s = \int_{\partial D} \abs{f_m}^2 \abs{\Phi'}^{-1} \, \di s < \infty 
\end{equation*}
since the singularities of $\abs{\Phi'}^{-1}$ at the preimages of the corners on $\partial \Omega$ are integrable~\cite[Theorem~3.2 \& Exercise~3.4.1]{Pommerenke92}. A similar argument also yields
\begin{equation*}
    \int_{\partial \Omega}  \widetilde{f}_m \, \di s =  \int_{\partial D} f_m\, \di s = 0.
\end{equation*}
To summarize, $\widetilde{f}_m \in L^2_\diamond(\partial \Omega)$ is a proper boundary current on $\partial \Omega$ for all $m \in \Z'$.

Assume that $\Omega$ is characterized by the conductivity $\widetilde{\gamma} \in L^\infty_+(\Omega)$ and drive the current density $\widetilde{f}_m$ through $\partial \Omega$. The resulting potential $\widetilde{u}_m^{\widetilde{\gamma}}$ in $\Omega$ is the unique weak solution of 
\begin{equation*}
    -\nabla \cdot (\widetilde{\gamma} \nabla \widetilde{u}) = 0 \quad \text{in } \Omega, \qquad  \nu \cdot (\widetilde{\gamma} \nabla \widetilde{u}) = \widetilde{f}_m \quad \text{on } \partial \Omega
\end{equation*}
in $H^1_\diamond(\Omega)$. According to \cite[Lemma~4.1 \& Remark~4.1]{Hyvonen17d}, the conformally mapped potential 
\begin{equation*}
    u^\gamma_m \coloneqq \widetilde{u}^{\widetilde{\gamma}}_m \circ \Phi
\end{equation*}
is a weak solution in $D$ to \eqref{eq:strong_form} for $f = f_m$ and $\gamma = \widetilde{\gamma} \circ \Phi$.

These observations can be interpreted as guidelines for simulating (or measuring) relative boundary data on $\Omega$ with the aim to reconstruct information on the conductivity perturbation $\widetilde{\eta} \coloneqq \widetilde{\gamma} - 1$:
\begin{enumerate}[(i)]
    \item Apply the current patterns $\{ \widetilde{f}_m \}_{m \in \Z_M'}$ from \eqref{eq:conformalfm} on $\partial \Omega$ for conductivities $\widetilde{\gamma}$ and $1$, and measure the resulting relative boundary potentials 
    \begin{equation*}
        \widetilde{g}_m^{\widetilde{\eta}} \coloneqq (\widetilde{u}^{\widetilde{\gamma}}_m - \widetilde{u}^{1}_m)|_{\partial \Omega}, \qquad m \in \Z_{M}'.
    \end{equation*}
    \item Form 
    \begin{equation*}
        g^\eta_m \coloneqq \widetilde{g}_m^{\widetilde{\eta}} \circ \Phi, \qquad m \in \Z_{M}',
    \end{equation*}
    that are relative boundary potentials on $\partial D$ for the Fourier current basis $\{ f_m \}_{m \in \Z_M'}$, corresponding to the conductivity perturbation $\eta = \widetilde{\gamma}\circ \Phi - 1 = \gamma - 1$ in $D$.
    \item Form the nonlinear data matrix
    \begin{equation} \label{eq:nonlindata}
        a_{m,n}^{\rm NL}(\eta) \coloneqq \inner{g^\eta_m, f_n}_{L^2(\partial D)}, \qquad m, n \in \Z_M'.
    \end{equation}
    \item Use the matrix $\{ a_{m,n}^{\rm NL}(\eta) \}_{m,n \in \Z_M'}$ in place of $\{ a_{m,n}(\eta) \}_{m,n \in \Z_M'}$ as the input for the algorithms of Section~\ref{sec:implementation} to reconstruct an approximation for~$\eta$.
    \item Form the final reconstruction of $\widetilde{\eta}$ in $\Omega$ by composing the reconstruction in $D$ with $\Psi$. 
\end{enumerate}

All employed CM data for the unit disk,~i.e.,~$\Omega = D$, are simulated by choosing $\Phi = {\rm Id}$ in the above scheme. Hence, also with $D$ as the examined domain, the data used as the input for the reconstruction algorithms of Section~\ref{sec:implementation} correspond to a nonlinear forward model in our numerical tests.
 
\subsection{CEM measurements}
\label{sec:CEM}

Assume next that the measurements of current and voltage on the boundary of the simply-connected domain $\Omega$ are performed with the help of $2L+1 \in \N \setminus \{ 1,2 \}$ contact electrodes $\{ E_l \}_{l \in \Z_L}$ that are modeled as open nonempty connected subsets of $\partial \Omega$ with mutually disjoint closures. That is, one can only control the mean-free net current pattern $I \in \C^{2L+1}$ through the electrodes and measure the resulting electrode potential vector $U \in \C^{2L +1}$. For simplicity, let us assume that all electrodes have the same arclength $d>0$. The thin resistive layers at the electrode object interfaces are modeled by contact impedances $\{ z_l \}_{l \in \Z_L} \subset \C$ that are assumed to have strictly positive real parts and to stay the same between measurements. The electrode measurements are modeled by the CEM~\cite{Cheng89,Somersalo92}, which is in good agreement with practical EIT.

As in Section~\ref{sec:piecewise_smooth}, suppose that we aim to reconstruct a conductivity perturbation $\widetilde{\eta} = \widetilde{\gamma} - 1 \in L^\infty(\Omega)$ by performing electrode measurements on $\partial \Omega$ with $\widetilde{\gamma}$ and $1$ acting in turns as the internal conductivity of $\Omega$. Assuming one is able to use a linearly independent set $\{ I_l \}_{l \in \Z_L'}$ of $2 L$ mean-free current patterns and freely choose the positions of the available electrodes, the construction presented in \cite[Section~6.2]{Garde21} produces an approximation $\widehat{g}_m^{\widetilde{\eta}}$ based on the measured electrode potentials for the relative boundary potential $\widetilde{g}_m^{\widetilde{\eta}}$ considered in Section~\ref{sec:piecewise_smooth}. If it is further assumed that $\partial \Omega$ is smooth and $\supp\widetilde{\eta}$ is compactly contained in $\Omega$, then \cite[Corollary~6.2]{Garde21} implies
\begin{align} \label{eq:CEM_approx}
    \norm{\widetilde{g}_m^{\widetilde{\eta}} - \widehat{g}_m^{\widetilde{\eta}}}_{L^2(\partial \Omega)} 
    &\leq C \bigl( \bigl (\tfrac{1}{2} + L\bigr)^{-s} + d^2\bigr) \norm{\widetilde{f}_m}_{H^s(\partial \Omega)}  \nonumber \\[1mm]
    &\leq C \bigl( \bigl (\tfrac{1}{2} + L\bigr)^{-s} + d^2\bigr) \norm{f_m}_{H^s(\partial D)} \nonumber \\[1mm]
    &\leq C \bigl( \bigl (\tfrac{1}{2} + L\bigr)^{-s} + d^2\bigr) \abs{m}^s
\end{align}
for any $s > 1/2$, with $C > 0$ being a generic constant (not depending on $m$) that may change between occurrences. The penultimate step follows from $\Psi|_{\partial \Omega}$ being smooth due to the smoothness of $\partial\Omega$~\cite{Pommerenke92}.

When the reconstruction algorithms of Section~\ref{sec:implementation} are applied to electrode measurements, we replace $\widetilde{g}_m^{\widetilde{\eta}}$ by $\widehat{g}_m^{\widetilde{\eta}}$ in the guidelines of Section~\ref{sec:piecewise_smooth}. According to \eqref{eq:CEM_approx}, the approximation provided by $\widehat{g}_m^{\widetilde{\eta}}$ is good for low Fourier frequencies $m$ and small electrodes. In particular, the parameter $M \in \N$ controlling the size of the Fourier data matrix should be chosen smaller than $L$ for \eqref{eq:CEM_approx} to guarantee a reasonable correspondence between $\widehat{g}_m^{\widetilde{\eta}}$ and $\widetilde{g}_m^{\widetilde{\eta}}$ for all involved Fourier frequencies $m \in \Z_M'$.

\begin{remark} \label{remark:CEM_meas}
    In our numerical experiments with electrode measurements, the imaged object is the unit disk with equiangular electrodes for simulated data and a water tank with the shape of a right circular cylinder for real-world measurements~\cite{Hauptmann17}. The considered conductivity distributions inside the water tank are homogeneous in the vertical direction, as are the rectangular electrodes that are equiangularly attached to lateral surface of the tank and extend from its bottom to the water surface. Due to certain symmetries, such three-dimensional water tank measurements can be modeled with the two-dimensional CEM in the unit disk, cf.,~e.g.,~\cite{Hyvonen17b}.

    Hence, one only needs to consider approximating the relative boundary potential $g_m^\eta = \widetilde{g}_m^{\widetilde{\eta}}$ from measurements on equiangular electrodes on the boundary of the unit disk. This is achieved by introducing a mean-free current pattern $I_m$ via evaluating a suitably scaled version of $f_m$ at the midpoints of the electrodes, measuring the resulting relative electrode potentials corresponding to the conductivity perturbation $\eta$, and finally performing trigonometric interpolation of the relative potentials with respect to the midpoints of the electrodes. See \cite[Sections~5 and 6.2]{Garde21} for the details.
\end{remark}
    
\section{Numerical experiments} \label{sec:res}

We employ three numerical techniques for simulating the relative potentials $g_m^\eta$, $m \in \Z_M'$ on the boundary of the unit disk $D$ corresponding to a given conformally mapped perturbation $\widetilde{\eta} \in L^\infty(\Omega)$, with $\widetilde{\eta} > -1$, in the imaged domain $\Omega$. Each of them generates approximate point values of $g_m$ on an equidistant grid with $2N+1$ nodes on $\partial D$. Subsequently, the nonlinear data matrix $\{ a_{m,n}^{\rm NL}(\eta) \}_{m,n \in \Z_M'}$ is computed via {\em fast Fourier transform} (FFT), which corresponds to approximating the integrals in \eqref{eq:nonlindata} by the periodic trapezoidal rule that has favorable convergence properties. The values of $M$ and $N$ vary between different tests.

If $\widetilde{\eta}= \eta$ is a scalar multiple of the characteristic function of a disk compactly contained in $\Omega = D$, the relative data $g_m^{\eta}$ is simulated by exploiting M\"obius transformations. This technique is based on $\{ f_m \}_{m \in \Z'}$ being the eigenfunctions with explicitly known eigenvalues for $\Lambda(1 + \kappa \chi_{D_{\rho}}) -  \Lambda(1)$, where $\kappa>-1$ is a constant and $\chi_{D_{\rho}}$ is the characteristic function of the open disk $D_\rho$ of radius $0 < \rho < 1$ centered at the origin. See \cite[Example~7.2]{Garde21} and \cite[Example~4]{Hyvonen18} for more details. This approach leads to highly accurate data, enabling studying the limits of the reconstruction method with noiseless data as well as the effect of the linearization error. Alternatively, one could use exact characterizations such as those in~\cite[Section~3]{Garde18}.

For other types of $\widetilde{\eta}$ and/or domains $\Omega$, the absolute potentials on $\partial \Omega$ for the current patterns $\{ \widetilde{f}_m \}_{m \in \Z_M'}$ corresponding to the conductivities $1 + \widetilde{\eta}$ and $1$ are simulated by a {\em finite element method} (FEM) with piecewise linear basis functions using the Python package scikit-fem \cite{skfem}. If $\Omega = D$, the {\em finite element} (FE) mesh is constructed so that the boundary nodes directly form an equidistant grid, which makes collecting the approximate nodal values of the relative potentials $g_m^{\eta}$ straightforward.  If $\Omega \not= D$ is a simply-connected polygon, the Schwarz--Christoffel toolbox \cite{Driscoll96} is utilized to numerically introduce the needed Riemann mapping $\Psi$. The boundary nodes of the FE mesh in $\Omega$ are then chosen as the images of an equidistant grid on $\partial D$ under $\Phi = \Psi^{-1}$. The simulated point values of the relative boundary potentials $\widetilde{g}^{\widetilde{\eta}}_m$, $m \in \Z_M'$, at the boundary nodes of the FE mesh are then precisely the values of $g_m^{\eta}$, $m \in \Z_M'$, on the equidistant grid of $\partial D$, cf.~Section~\ref{sec:piecewise_smooth}.

When considering electrode measurements on $\Omega = D$, the CEM forward problems for the conductivities $1+ \eta$ and $1$ are numerically solved with FEM based on piecewise linear basis functions for the potential in $\Omega$ and certain macro elements for the electrode voltages~\cite{Vauhkonen97}. Constructing an approximation $\widehat{g}_m^{\eta}$ for the relative boundary potential $\widetilde{g}_m^{\widetilde{\eta}} = g_m^{\eta}$ is performed as explained in Remark~\ref{remark:CEM_meas}, with the midpoints of the electrodes serving as the uniform grid on $\partial D$, resulting in $N = L$. Note that only accounting for the point values of $\widehat{g}_m^{\eta}$ at the electrode centers when evaluating the integrals in~\eqref{eq:nonlindata}  does not lead to a reduction in accuracy, since the periodic trapezoidal rule is exact when computing~\eqref{eq:nonlindata} with $g_m^{\eta}$ replaced by $\widehat{g}_m^{\eta}$; see~\cite{Saranen02} and \cite[Eq.~(4.8)]{Garde21}.

When considering noisy simulated measurements, we add independent realizations of Gaussian noise to the real and imaginary parts of all elements in the nonlinear data matrix $\{ a_{m,n}^{\rm NL}(\eta) \}_{m,n \in \Z_M'}$. We say that the data contains $100\sigma\%$ of additive noise, with $\sigma\in[0,1]$, if the noise added to $\redel(a_{m,n}^{\rm NL})$ is drawn from $\mathcal{N}(0, \sigma^2 (\redel(a_{m,n}))^2  )$ and that added to $\imdel(a_{m,n}^{\rm NL})$ is drawn from $\mathcal{N}(0, \sigma^2 (\imdel(a_{m,n}))^2  )$ for all $m,n \in \Z_M'$. For simulated measurements, we assume to know the statistics of the additive noise and choose $\delta$ in \eqref{eq:morozov1} and \eqref{eq:morozov2} to be the square root of the expected value for the squared Euclidean
norm of the noise vector, i.e., the square root of the sum of the variances of the independent noise components. On the other hand, the value for the fudge factor $\omega$ in \eqref{eq:morozov1} and \eqref{eq:morozov2} is case-dependent.

\subsection{Visualization of right singular vectors of $F^M$}

We start the numerical experiments by computing the right singular vectors of $F^M$ for $M=32$, interpreting them as coefficients in a truncated Zernike polynomial expansion via \eqref{eq:finite_complete_system}, \eqref{eq:finite_triang_syst}, \eqref{eq:inf_vec} and \eqref{eq:Zexpansion}, and visualizing the resulting right singular functions over $D$. We do not prove that the singular functions or values of $F^M$ approximate such for the untruncated operator $F$. We also remark that the order in which singular functions with certain qualitative properties appear in the SVD depends on $M$. 

Due to the block diagonal structure \eqref{eq:bloc_triangular}, each right singular function of $F^M$ is a linear combination of Zernike polynomials with a fixed absolute angular index $\abs{j}$. In fact, there are $M$ singular functions for $j=0$ and $2(M-\abs{j})$ for each $\abs{j} \in \{1,\dots,M-1\}$. The first right singular function is a linear combination of Zernike polynomials for $j=0$ over different radial indices $k$. The second largest singular value has multiplicity two, and the corresponding two orthonormal singular functions can be chosen to be linear combinations of Zernike polynomials for $j=1$ and $j=-1$, respectively, with the same coefficients for each $k$ in \eqref{eq:Zexpansion}. This is the convention we choose for any fixed $\abs{j}$ in what follows. The next $13$ singular values also have multiplicity two, and they correspond to $\abs{j} = 2, \dots, 14$. Tables~\ref{tbl:sv_real} and \ref{tbl:sv_imag}, respectively, show the real and imaginary parts of the first nine right singular functions that exhibit higher angular frequencies and faster decay toward the origin as $\abs{j}$ increases. 

\begin{table}[htb]
    \centering
    \begin{tabular}{ccccc}
        $j=0$ & $j=\pm 1$ & $j=\pm 2$ & $j=\pm 3$ & $j=\pm 4$ \\
        \midrule 
         \includegraphics[width=0.15\textwidth]{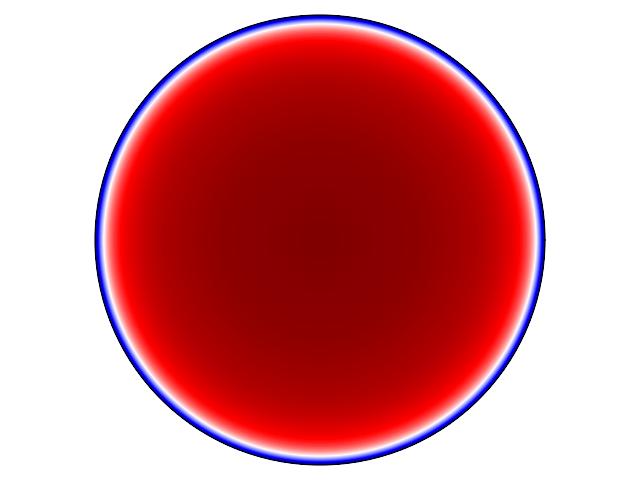}
         & 
        \includegraphics[width=0.15\textwidth]{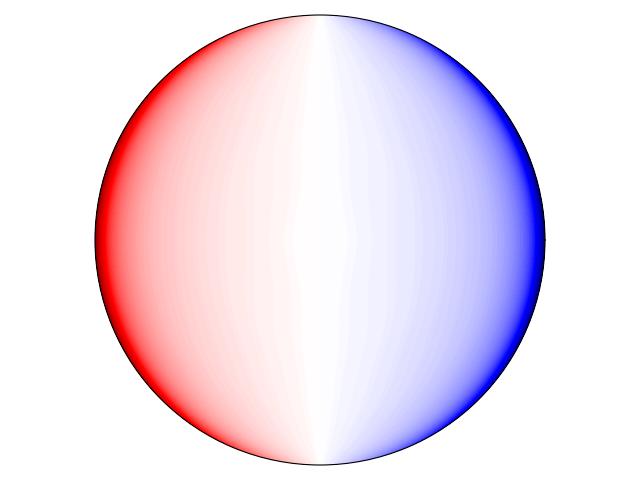}
        & 
        \includegraphics[width=0.15\textwidth]{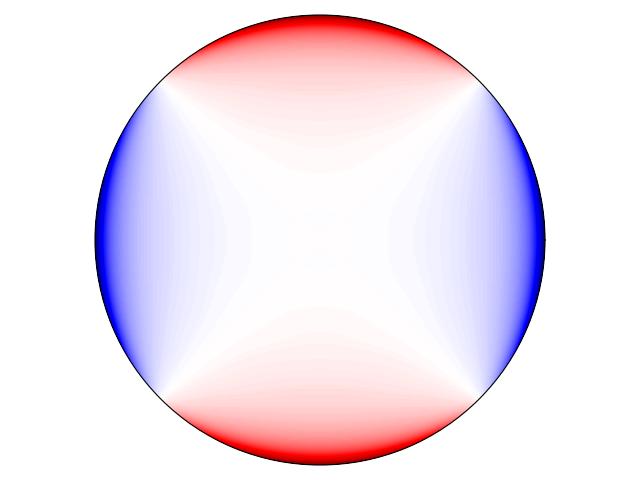}
        & 
        \includegraphics[width=0.15\textwidth]{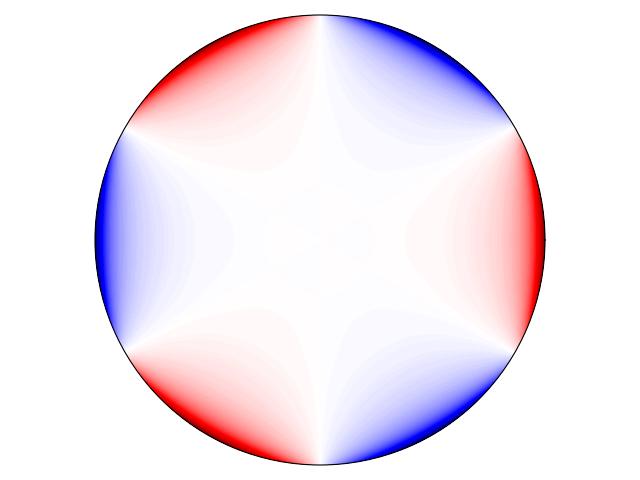}
        &
        \includegraphics[width=0.15\textwidth]{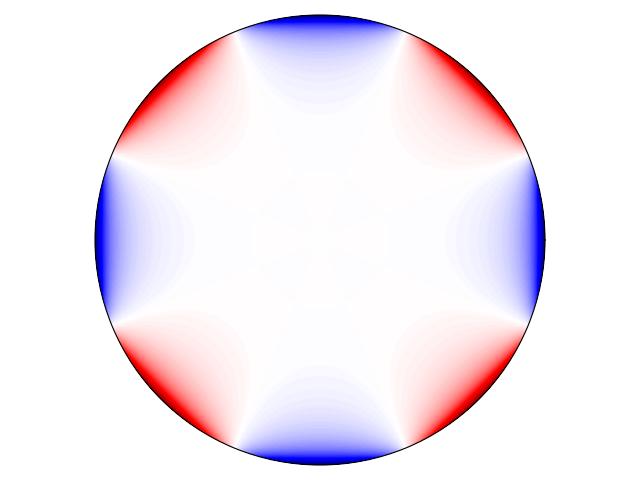}
        \\
         & 
        \includegraphics[width=0.15\textwidth]{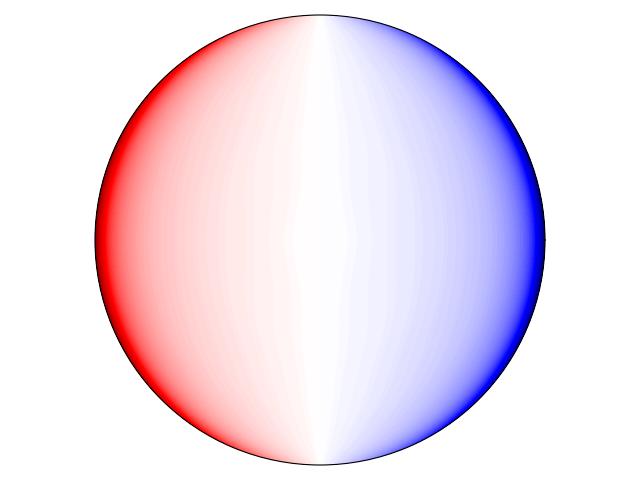}
        & 
        \includegraphics[width=0.15\textwidth]{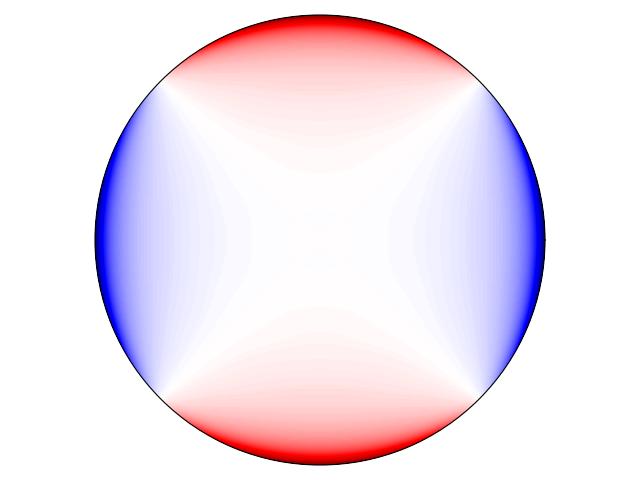}
        & 
        \includegraphics[width=0.15\textwidth]{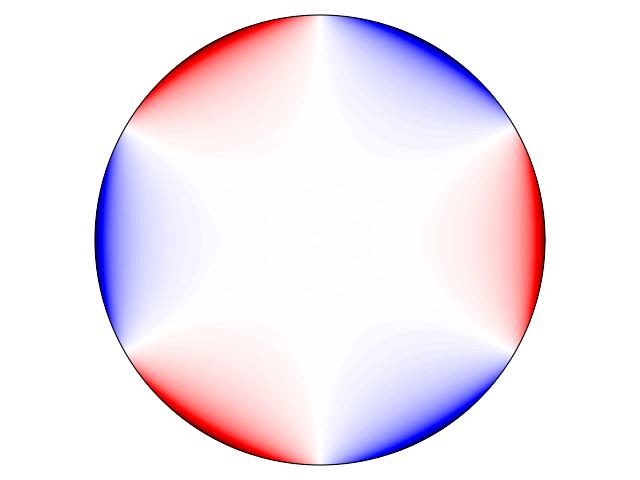}
        &
        \includegraphics[width=0.15\textwidth]{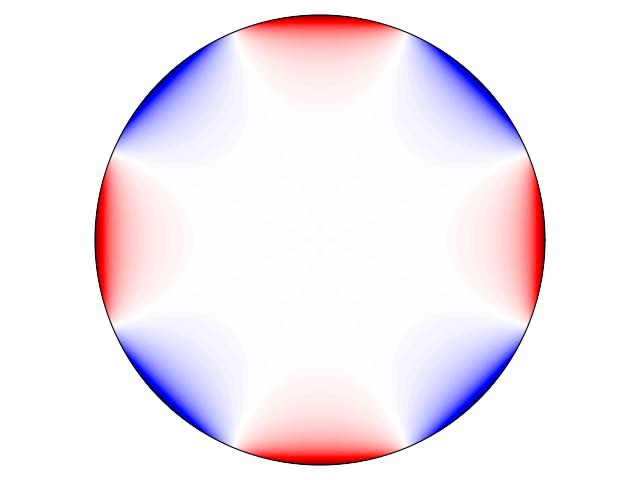}
    \end{tabular}
    \caption{The real parts of the first nine right singular functions of $F^M$ for $M=32$.}
    \label{tbl:sv_real}
\end{table}

\begin{table}[htb]
    \centering
    \begin{tabular}{ccccc}
        $j=0$ & $j=\pm 1$ & $j=\pm 2$ & $j=\pm 3$ & $j=\pm 4$ \\
        \midrule 
         \includegraphics[width=0.15\textwidth]{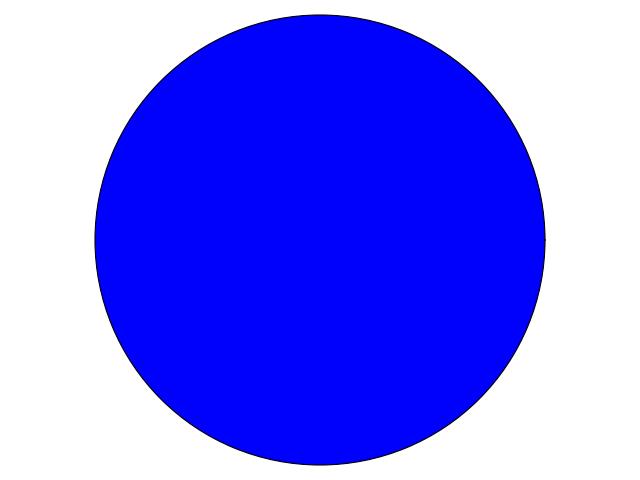}
         & 
        \includegraphics[width=0.15\textwidth]{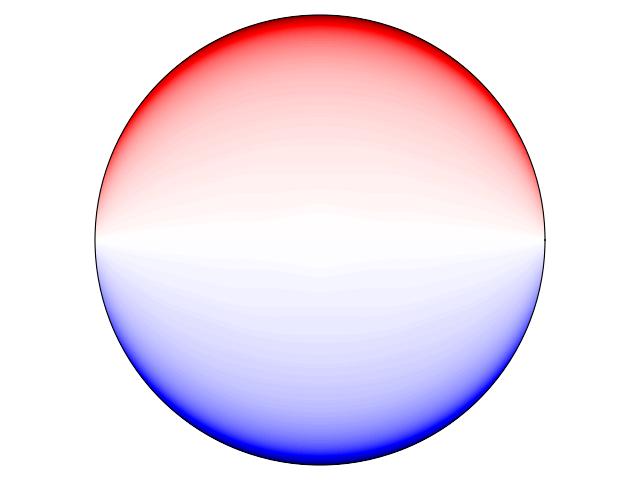}
        & 
        \includegraphics[width=0.15\textwidth]{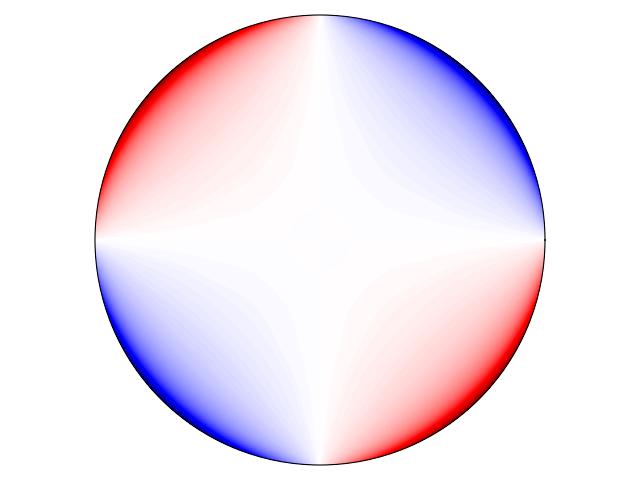}
        & 
        \includegraphics[width=0.15\textwidth]{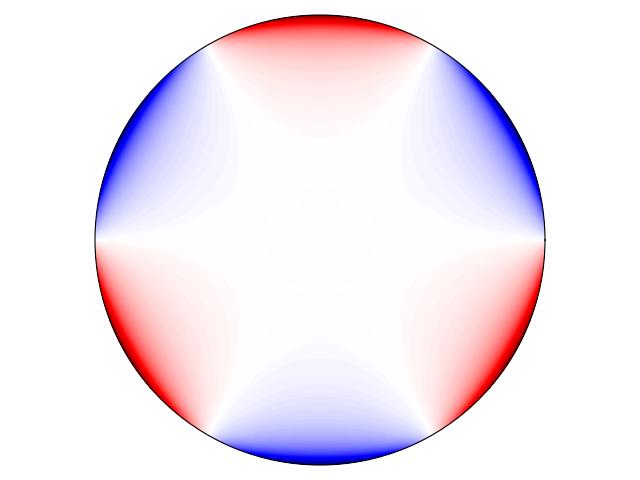}
        &
        \includegraphics[width=0.15\textwidth]{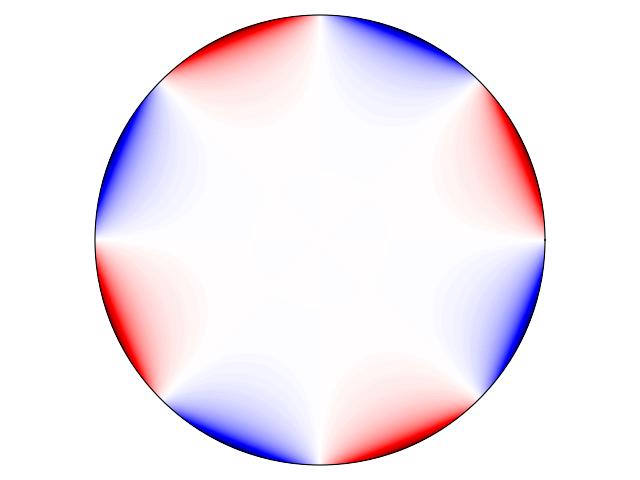}
        \\
         & 
        \includegraphics[width=0.15\textwidth]{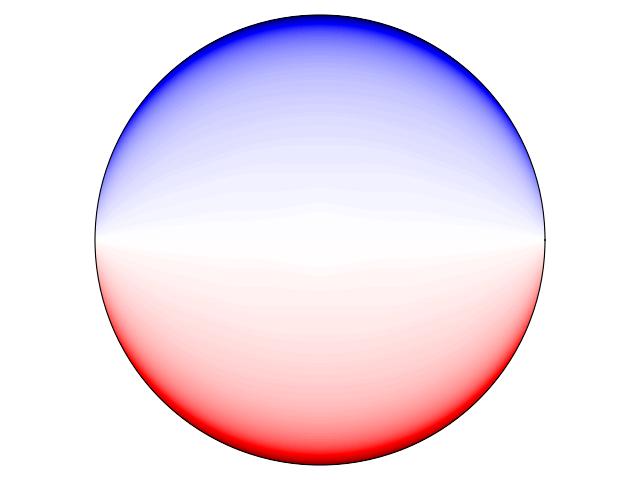}
        & 
        \includegraphics[width=0.15\textwidth]{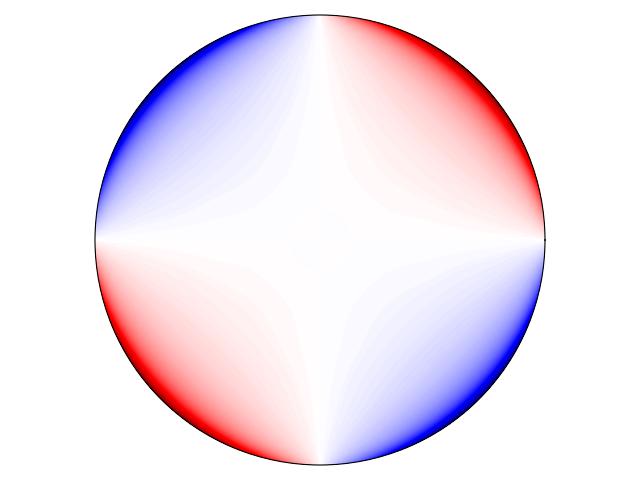}
        & 
        \includegraphics[width=0.15\textwidth]{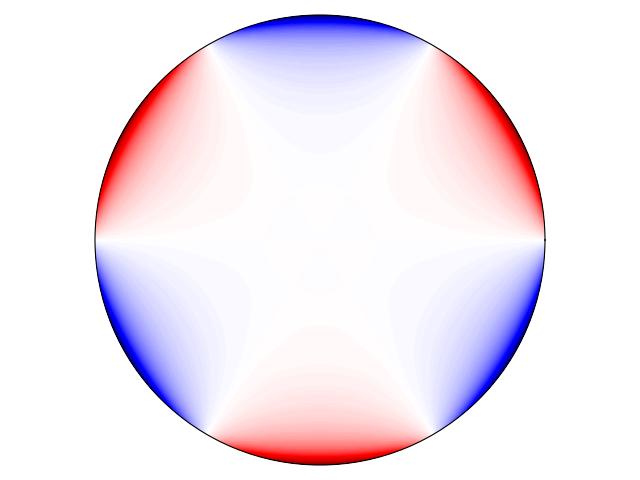}
        &
        \includegraphics[width=0.15\textwidth]{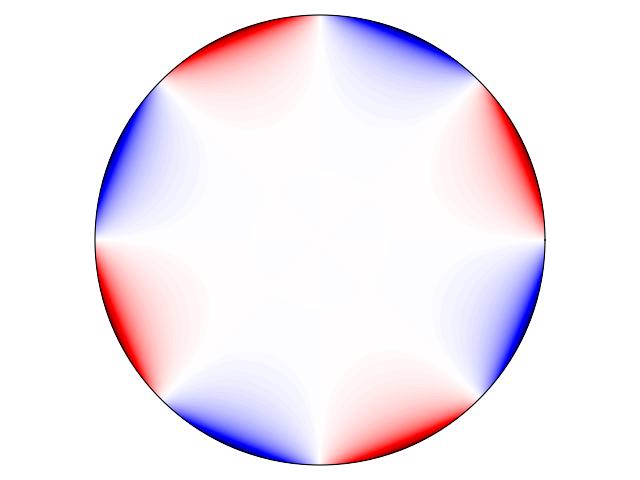}
    \end{tabular}
    \caption{The imaginary parts of the first nine right singular functions of $F^M$ for $M=32$.}
    \label{tbl:sv_imag}
\end{table}

At a certain point, singular functions that are new linear combinations of Zernike polynomials over $k$ for an already visited $\abs{j}$ start to appear. For our choice of $M=32$, the first such singular function is the 30th that corresponds to the lowest angular frequency $j=0$. The next ones are the 45th and 46th, 53rd and 54th, as well as 59th and 60th singular functions that are associated to $j=\pm 1$, $j=\pm 2$ and $j=\pm 3$, respectively. Their real and imaginary parts, illustrated in Tables~\ref{tbl:sv_real_more} and \ref{tbl:sv_imag_more}, exhibit similar angular oscillations as the first singular functions for the same $\abs{j}$ in Tables~\ref{tbl:sv_real} and \ref{tbl:sv_imag}. However, these `second order' singular functions in Tables~\ref{tbl:sv_real_more} and \ref{tbl:sv_imag_more} carry an additional radial oscillation close to the boundary of $D$.

\begin{table}[htb]
    \centering
    \begin{tabular}{cccc}
        30th & 45th and 46th & 53rd and 54th & 59th and 60th \\
        singular function,& singular functions,& singular functions,& singular functions,\\
        $j=0$ & $j=\pm 1$ & $j=\pm 2$ & $j=\pm 3$\\
        \midrule 
         \includegraphics[width=0.15\textwidth]{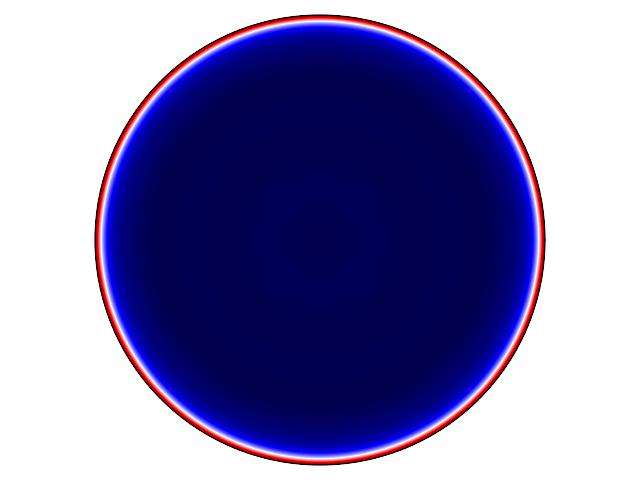}
         & 
        \includegraphics[width=0.15\textwidth]{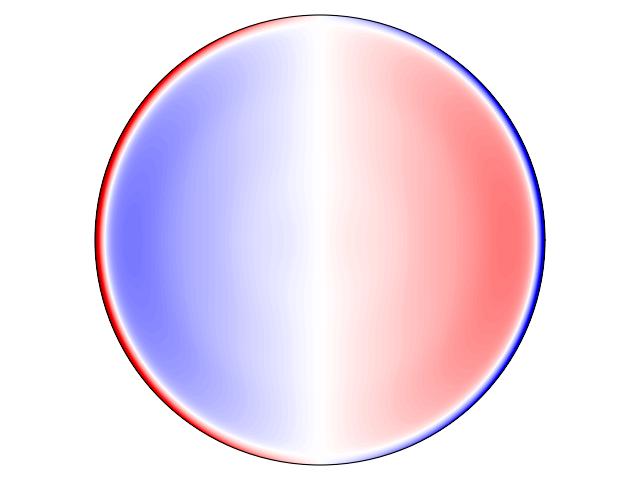}
        & 
        \includegraphics[width=0.15\textwidth]{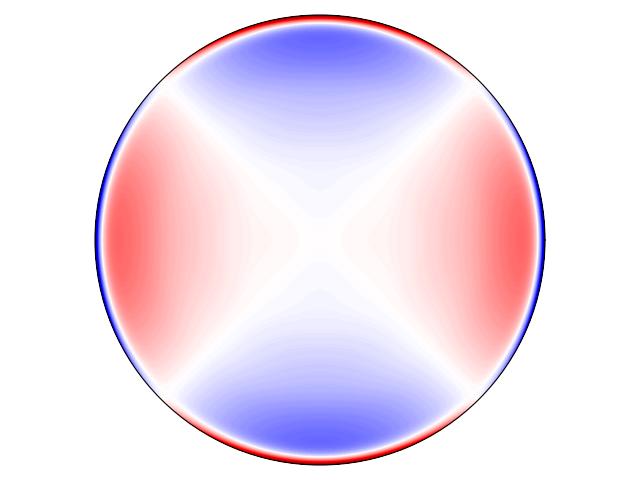}
        & 
        \includegraphics[width=0.15\textwidth]{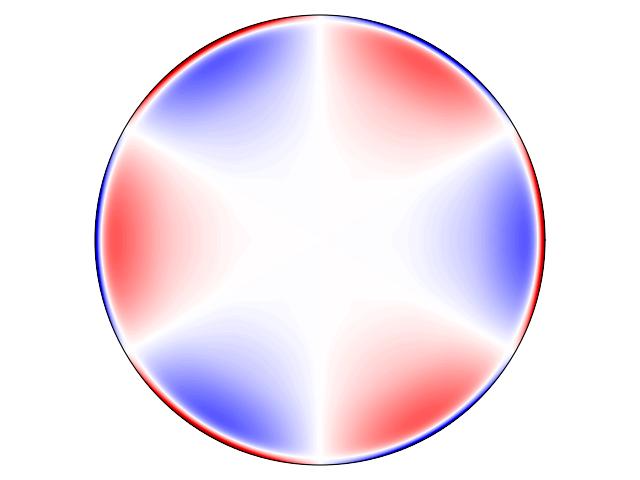}
        \\
        & 
        \includegraphics[width=0.15\textwidth]{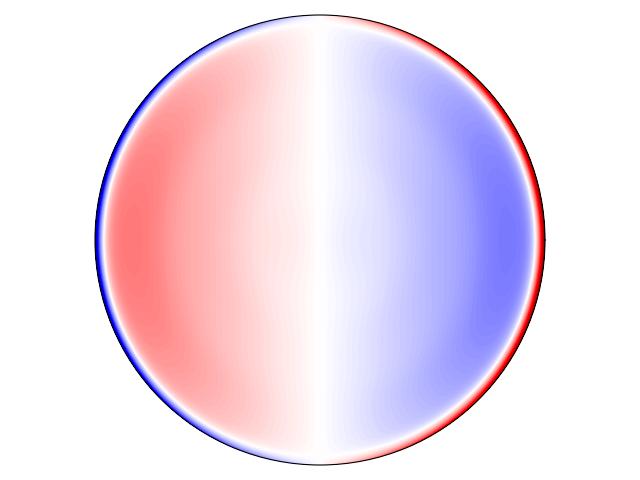}
        & 
        \includegraphics[width=0.15\textwidth]{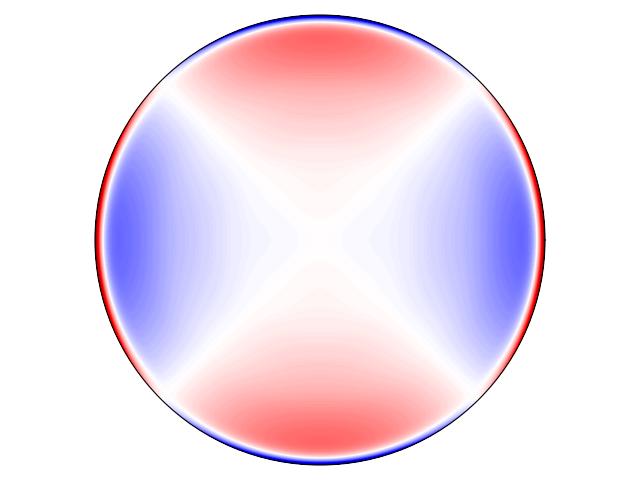}
        &
        \includegraphics[width=0.15\textwidth]{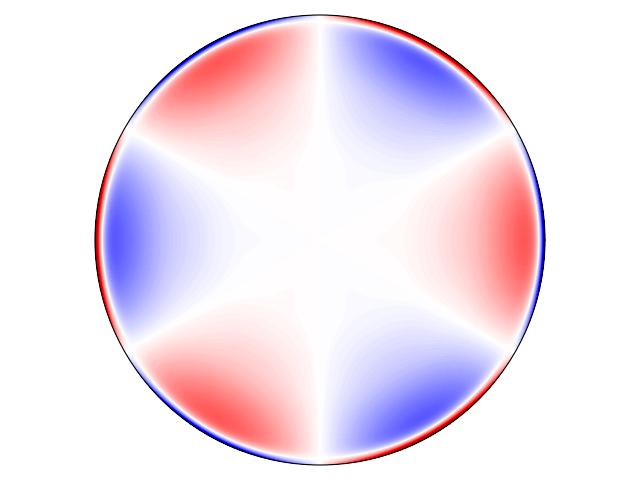}
    \end{tabular}
    \caption{The real parts of the first seven right singular functions of $F^M$ that are linear combinations of Zernike polynomials for some already visited angular index $j$ with $M=32$.}
    \label{tbl:sv_real_more}
\end{table}

\begin{table}[htb]
    \centering
    \begin{tabular}{cccc}
        30th & 45th and 46th & 53rd and 54th & 59th and 60th \\
        singular function,& singular functions,& singular functions,& singular functions,\\
        $j=0$ & $j=\pm 1$ & $j=\pm 2$ & $j=\pm 3$\\
        \midrule 
         \includegraphics[width=0.15\textwidth]{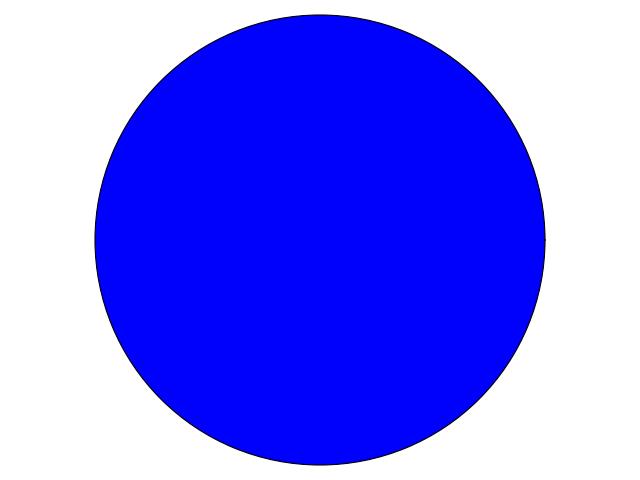}
         & 
        \includegraphics[width=0.15\textwidth]{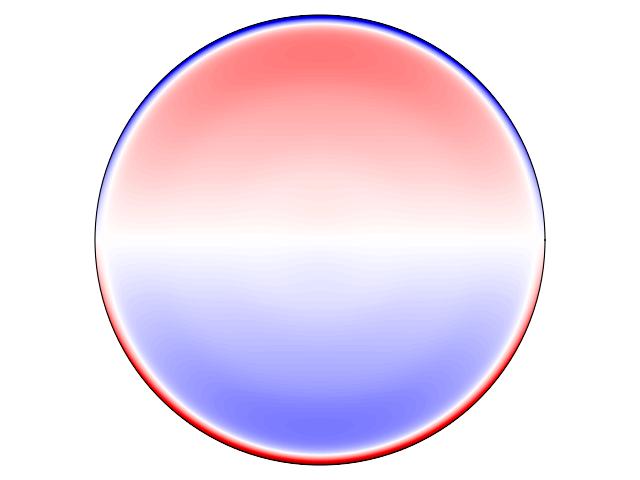}
        & 
        \includegraphics[width=0.15\textwidth]{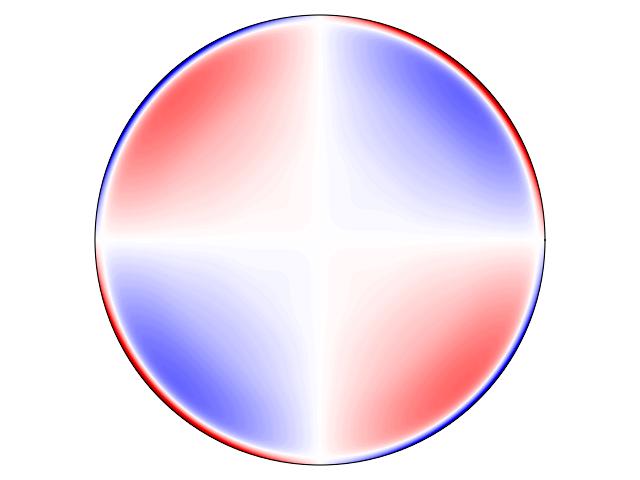}
        & 
        \includegraphics[width=0.15\textwidth]{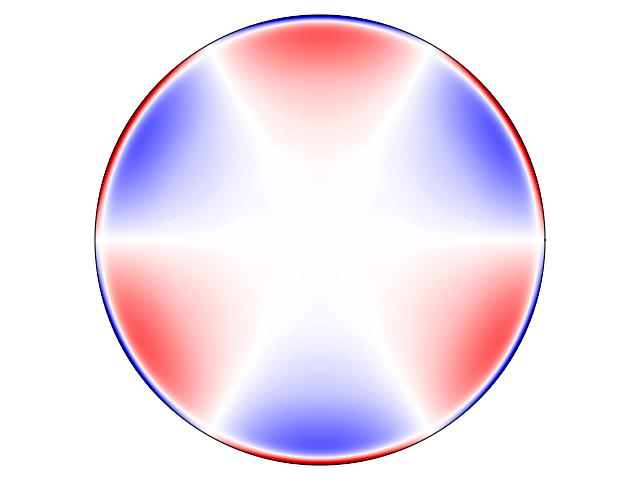}
        \\
        & 
        \includegraphics[width=0.15\textwidth]{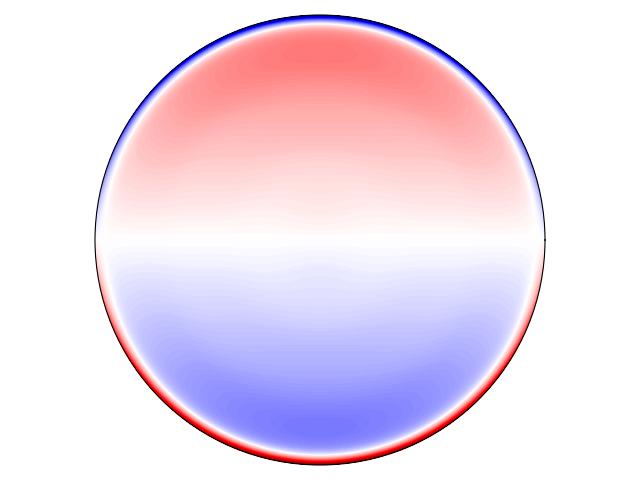}
        & 
        \includegraphics[width=0.15\textwidth]{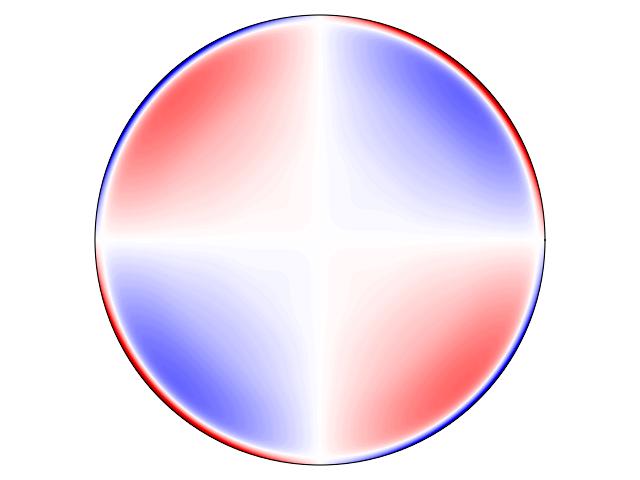}
        &
        \includegraphics[width=0.15\textwidth]{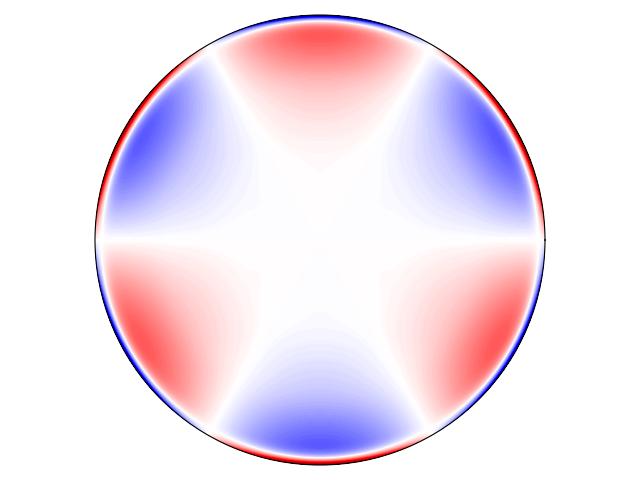}
    \end{tabular}
    \caption{The imaginary parts of the first seven right singular functions of $F^M$ that are linear combinations of Zernike polynomials for some already visited angular index $j$ with $M=32$.}
    \label{tbl:sv_imag_more}
\end{table}

\begin{figure}[h!]
     \centering
     \begin{subfigure}[t]{0.4\textwidth}
         \centering
         \includegraphics[width=\textwidth]{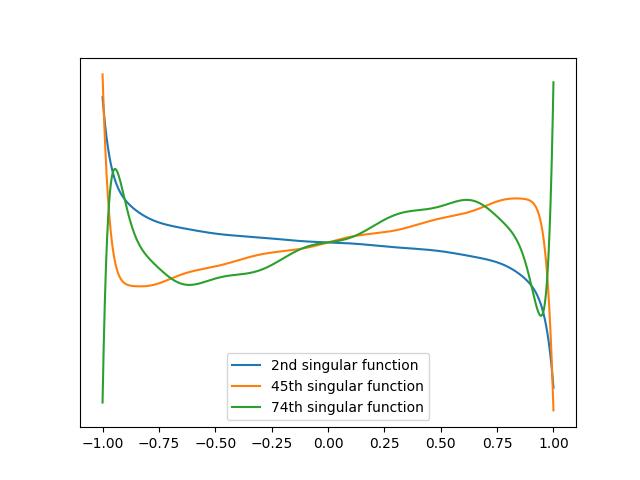}
         \subcaption{The real parts of the first three right singular functions for $j=1$ along the horizontal axis.}
     \end{subfigure}
     \hspace{2cm}
     \begin{subfigure}[t]{0.4\textwidth}
         \centering
         \includegraphics[width=\textwidth]{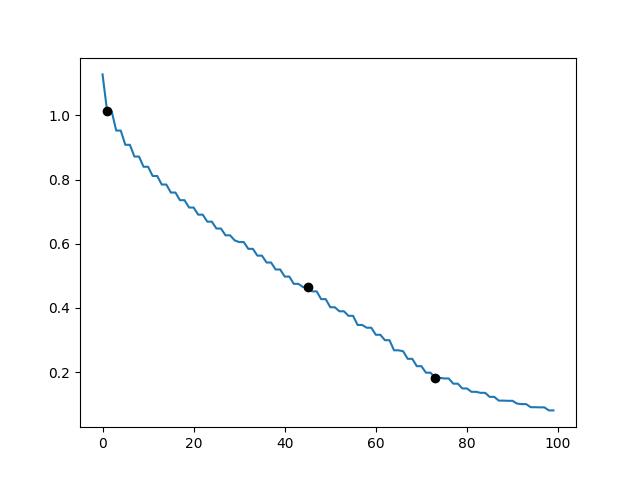}
         \subcaption{The first 100 singular values with the ones corresponding to the singular functions in the left-hand image highlighted.}
     \end{subfigure}
     \hfill
     \caption{The radial components of certain right singular functions and the first hundred singular values of $F^M$ for $M=32$.}
     \label{fig:sv_slice}
\end{figure}

As expected, the first `third order' singular function corresponds to $j=0$, and it is altogether the 67th one.  The next angular frequencies $j=\pm 1$ are visited for the third time at the 74th and 75th singular functions. As a general rule, the smaller the singular value corresponding to a given absolute angular index $\abs{j}$ is, the more oscillations there are in the radial direction in the associated singular functions. To illustrate this, the left-hand image of Fig.~\ref{fig:sv_slice} plots the real parts of the first three singular functions (2nd, 45th, 74th) for $j=1$ along the real axis. Finally, the right-hand image of Fig.~\ref{fig:sv_slice} presents, for the sake of completeness, the first 100 singular values of $F^M$ counted according to their multiplicity.

\subsection{CM and a discoidal inclusion in $D$} \label{sec:disk_in_disk}

\begin{figure}[htb]
     \centering
     \begin{subfigure}[t]{0.33\textwidth}
         \centering
         \includegraphics[width=\textwidth]{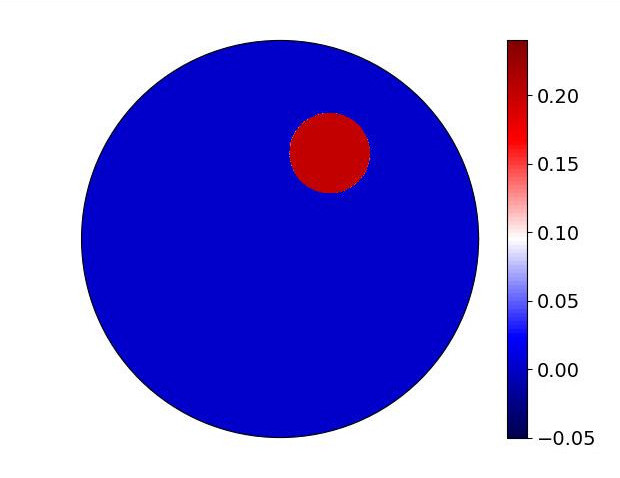}
         \subcaption{The target conductivity perturbation.}
         \label{fig:reco_disk_target}
     \end{subfigure}
     \qquad \qquad
     \begin{subfigure}[t]{0.33\textwidth}
         \centering
         \includegraphics[width=\textwidth]{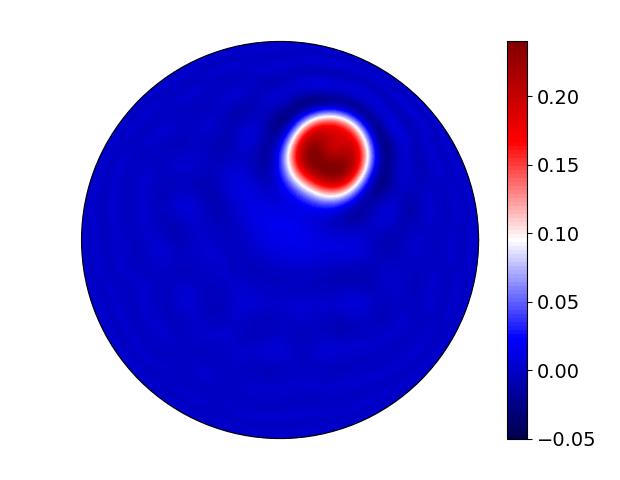}
         \subcaption{Very accurate data and minimal SVD-based regularization with $p=492$.}
         \label{fig:reco_disk_accu}
     \end{subfigure}
     \\
     \begin{subfigure}[t]{0.33\textwidth}
         \centering
         \includegraphics[width=\textwidth]{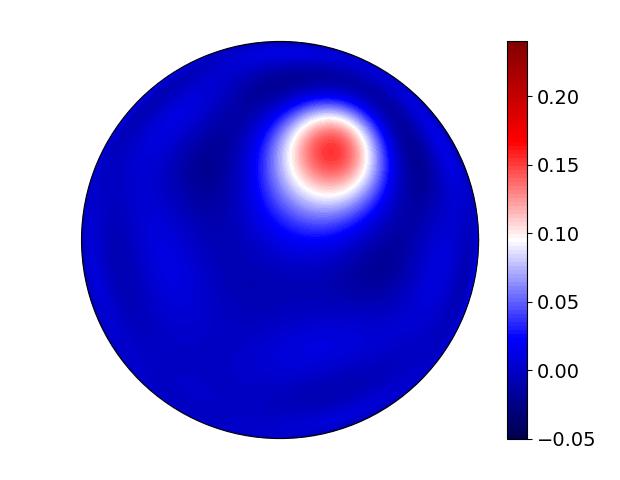}
         \subcaption{$1\%$ of  noise, SVD-based regularization with $\omega=1$, $p=174$.}
         \label{fig:reco_disk_svd}
     \end{subfigure}
    \qquad \qquad
     \begin{subfigure}[t]{0.33\textwidth}
         \centering
         \includegraphics[width=\textwidth]{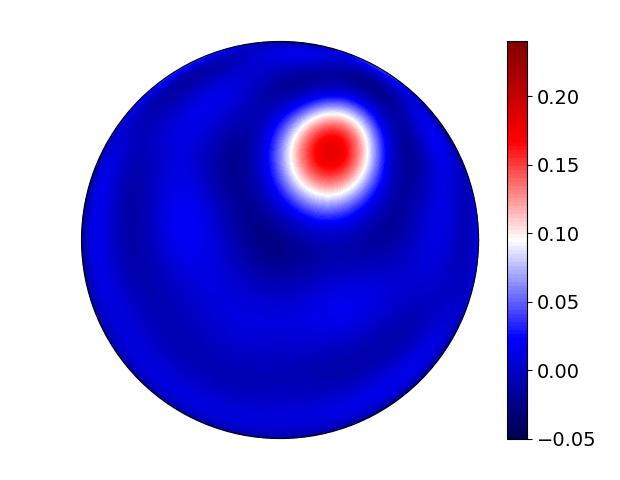}
         \subcaption{$1\%$ of noise, truncation of triangular subsystems with $\omega=1$, $p=154$.}
         \label{fig:reco_disk_hrc}
     \end{subfigure}
     \caption{Reconstructions of a discoidal inclusion inside the unit disk from highly accurate data generated using a M\"obius transformation with $M=32$.}
     \label{fig:reco_disk}
\end{figure}

The first numerical tests on the actual reconstruction algorithms consider the target in Fig.~\ref{fig:reco_disk_target}, corresponding to a discoidal inclusion of radius $0.2$ centered at $(\tfrac{1}{4},\tfrac{\sqrt{3}}{4})$ inside the unit disk. The conductivity of the inclusion is $1.2$ and, as in all our tests based on simulated measurements, that of the background is $1$. We visualize all targets and reconstructions as perturbations of the background conductivity, which explains the colormap of Fig.~\ref{fig:reco_disk}. Highly accurate data $\{ a_{m,n}^{\rm NL}(\eta) \}_{m,n \in \Z_M'}$  is generated for $M=32$ using M\"obius transformations as briefly explained at the beginning of Section~\ref{sec:res}. 

First, the truncated SVD algorithm described in Section~\ref{sec:svd_regularization} is run on the noiseless data with minimal truncation. According to a subjective visual inspection over many truncation indices, $p=492$ gives a close to optimal reconstruction that is shown in Fig.~\ref{fig:reco_disk_accu} and demonstrates that the support of a target inclusion can indeed be faithfully reconstructed despite the employed linearization if the available data is very accurate,~cf.~\cite{Harrach12}. This conclusion remains valid even for inclusions with higher contrasts. Although the reconstruction in Fig.~\ref{fig:reco_disk_accu} also reproduces the conductivity level of the inclusion rather accurately, this would no longer be the case if the conductivity of the inhomogeneity differed considerably more from the background value, due to the linearization error.

Next both the truncated SVD method and the algorithm based on truncation of triangular subsystems from Section~\ref{sec:triangular_regularization} are applied to data with $1\%$ of additive noise. We set $\omega = 1$ in \eqref{eq:morozov1} and \eqref{eq:morozov2}, which means that there is no attempt to compensate for the employed linearization or numerical errors via the choice of the fudge factor. The resulting reconstructions are shown in Fig.~\ref{fig:reco_disk_svd} and Fig.~\ref{fig:reco_disk_hrc}, and they correspond to the truncation indices $p=174$ and $p=154$ in the respective algorithms. In particular, far fewer Zernike coefficients are considered in Figs.~\ref{fig:reco_disk_svd} and \ref{fig:reco_disk_hrc} than in Fig.~\ref{fig:reco_disk_accu}. The two regularization methods seem to work about equally well, with the SVD-based approach better recapturing the unperturbed background but being worse at estimating the strength of the conductivity perturbation.

\begin{figure}[htb]
     \centering
     \begin{subfigure}[t]{0.3\textwidth}
         \centering
         \includegraphics[width=\textwidth]{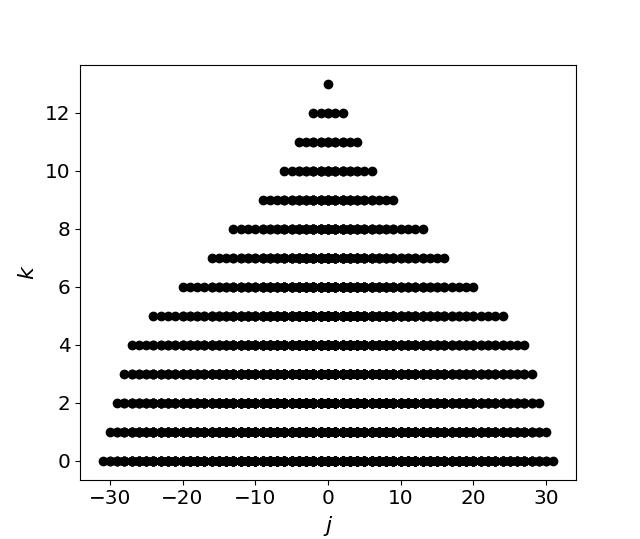}
         \subcaption{The SVD-based reconstruction in Fig~\ref{fig:reco_disk_accu}.}
     \end{subfigure}
     \hfill
     \begin{subfigure}[t]{0.3\textwidth}
         \centering
         \includegraphics[width=\textwidth]{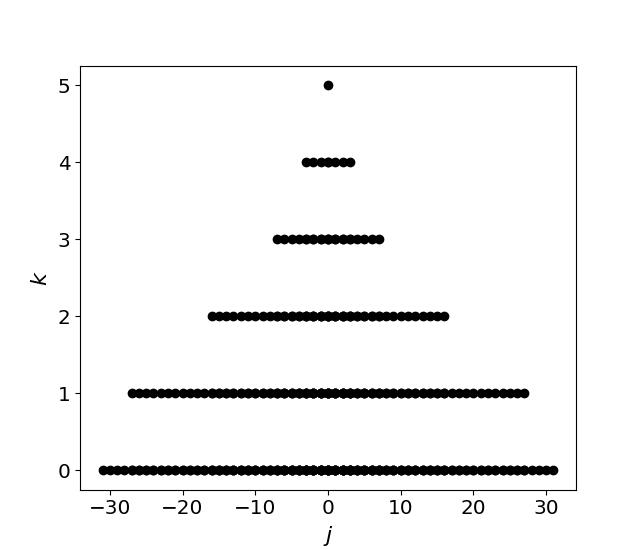}
         \subcaption{The SVD-based reconstruction in Fig~\ref{fig:reco_disk_svd}.}
     \end{subfigure}
     \hfill
     \begin{subfigure}[t]{0.3\textwidth}
         \centering
         \includegraphics[width=\textwidth]{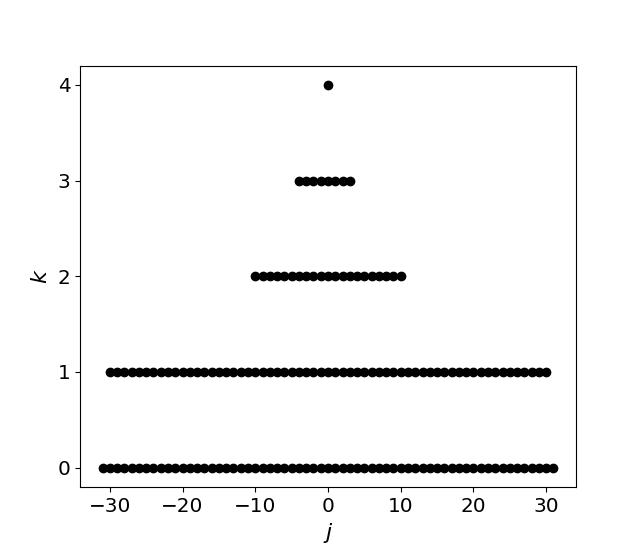}
         \subcaption{Truncation of triangular subsystems in  Fig~\ref{fig:reco_disk_hrc}.}
     \end{subfigure}
     \hfill
     \caption{The accepted Zernike indices for the reconstructions in Fig~\ref{fig:reco_disk}.}
     \label{fig:zernike_indices}
\end{figure}

The Zernike coefficients $c_{j,k}$ that contribute to the reconstructions in Figs.~\ref{fig:reco_disk_accu}, ~\ref{fig:reco_disk_svd} and ~\ref{fig:reco_disk_hrc}, i.e.,~the coefficients that are accepted/estimated by the corresponding regularization algorithms, are plotted in Fig.~\ref{fig:zernike_indices}. For both algorithms, the number of considered radial indices $k$ decreases as a function of $|j|$, with this effect being more pronounced for the truncated SVD method.

\subsection{CM for other types of perturbations and domains}

In this section the `nonlinear' data matrix $\{ a_{m,n}^{\rm NL}(\eta) \}_{m,n \in \Z_M'}$, with $M=16$, is simulated by FEM using $1024$ boundary nodes. This leads to FE meshes with approximately 130,000 nodes and 260,000 triangles for all considered geometries. 

\begin{figure}[b!]
     \centering
     \begin{subfigure}[t]{0.3\textwidth}
         \centering
         \includegraphics[width=\textwidth]{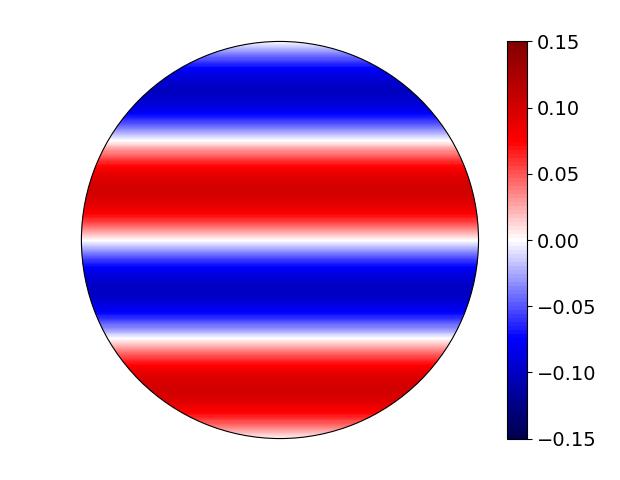}
         \subcaption{The target conductivity perturbation.}
          \label{fig:reco_sin_target}
     \end{subfigure}
     \hfill
     \begin{subfigure}[t]{0.3\textwidth}
         \centering         \includegraphics[width=\textwidth]{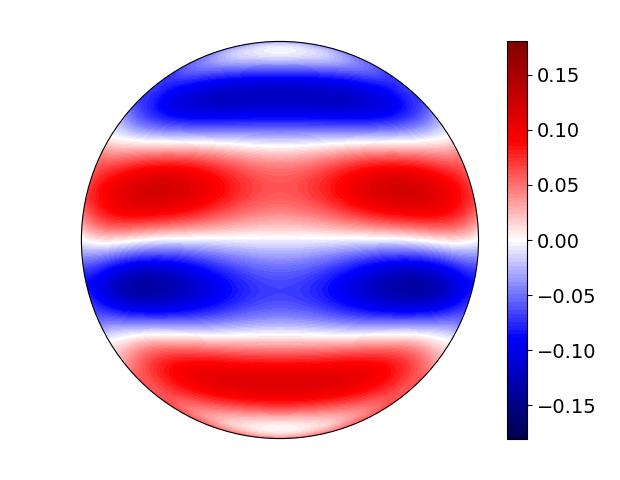}
         \subcaption{$1\%$ of noise, truncated SVD with $\omega=1$, $p=151$.
         }
     \end{subfigure}
     \hfill
     \begin{subfigure}[t]{0.3\textwidth}
         \centering
         \includegraphics[width=\textwidth]{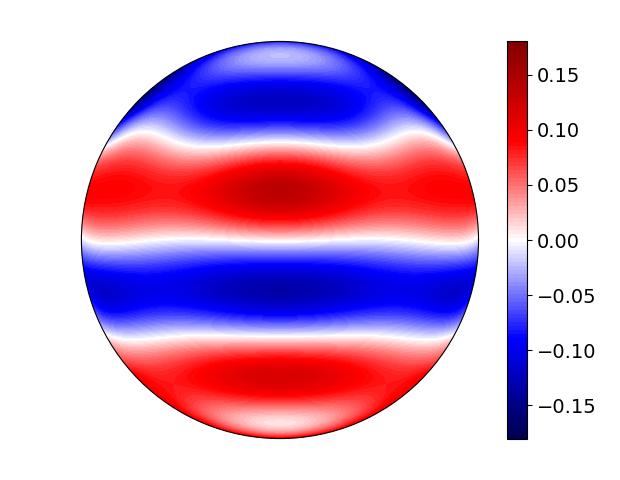}
         \subcaption{
         $1\%$ of noise, truncated triangular subsystems with $\omega=60$, $p=120$.}
     \end{subfigure}
     \hfill
     \\
     \begin{subfigure}[t]{0.3\textwidth}
         \centering
     \includegraphics[width=\textwidth]{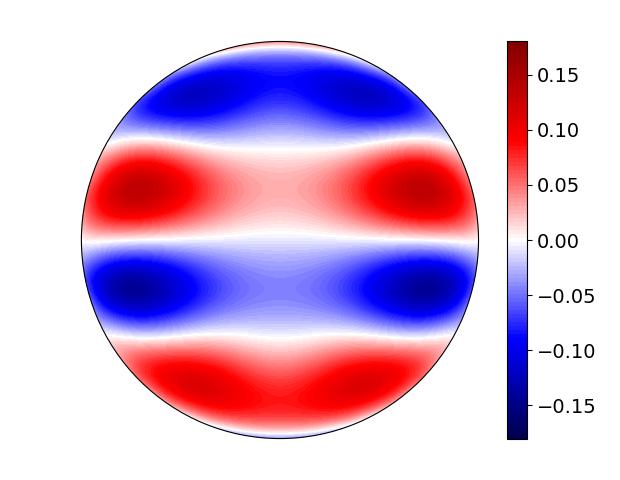}
         \subcaption{
         $10\%$ of noise, truncated SVD with $\omega=1$, $p=114$.
          }
     \end{subfigure}
     \qquad 
     \begin{subfigure}[t]{0.3\textwidth}
         \centering
    \includegraphics[width=\textwidth]{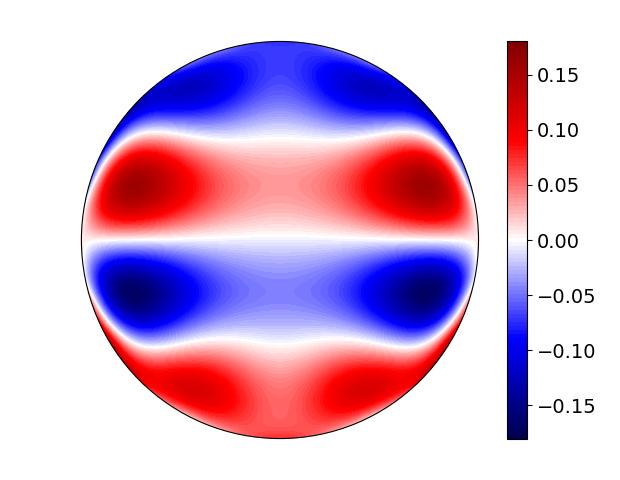}
         \subcaption{
         $10\%$ of noise, truncated triangular subsystems with $\omega=7$, $p=78$.
          }
     \end{subfigure}
    
     \caption{Reconstructions of a wave-like conductivity perturbation in the unit disk with $M=16$.}
     \label{fig:reco_sin}
\end{figure}

Figure \ref{fig:reco_sin_target} illustrates the unit disk with a target conductivity perturbation that behaves as a sine wave with amplitude 0.1 in the vertical direction and is homogeneous in the horizontal direction. The other four images in Fig.~\ref{fig:reco_sin} show reconstructions produced by the truncated SVD method (Section~\ref{sec:svd_regularization}) and truncation of triangular subsystems (Section~\ref{sec:triangular_regularization}) for the noise levels of 1\% and 10\%. The SVD-based reconstructions still correspond to the trivial choice $\omega = 1$ in~\eqref{eq:morozov1}. However, for the truncation of triangular subsystems, the fudge factor needs to be manually tuned in~\eqref{eq:morozov2} to enable reasonable reconstructions, leading to the heuristic choices $\omega = 60$ and $\omega = 7$, respectively, for the lower and higher noise levels. 

The reconstructions in Fig.~\ref{fig:reco_sin} are able capture the qualitative structure of the target conductivity for both noise levels, although the periodic pattern is somewhat distorted in all reconstructions. This effect is most visible close to the boundary and center of the domain, with reconstructions with the higher noise level suffering from stronger artifacts in this regard. Be that as it may, Fig.~\ref{fig:reco_sin} demonstrates that ignoring the linearization error may still lead to useful reconstructions in EIT even if the supports of the to-be-reconstructed conductivity perturbations are not restricted to the interior of the imaged domain.

A possible explanation for the requirement to tune the Morozov principle for the algorithm based on the truncation of triangular subsystems is its sensitivity to numerical inaccuracies in those elements of the data matrix that correspond to high spatial frequencies: in Section~\ref{sec:disk_in_disk} all simulated data were highly accurate, whereas the FEM simulations of this section are bound to suffer from numerical inaccuracies that are presumably the largest for the measurements corresponding to the highest spatial frequencies in the set of employed Fourier current patterns $\{ f_m \}_{m \in \Z_{16}'}$. As all our remaining tests consider either data simulated by FEM or/and CEM measurements, we only employ the SVD-based regularization method in what follows.

\begin{figure}[htb]
     \centering
     \begin{subfigure}[t]{0.4\textwidth}
         \centering
         \includegraphics[width=\textwidth]{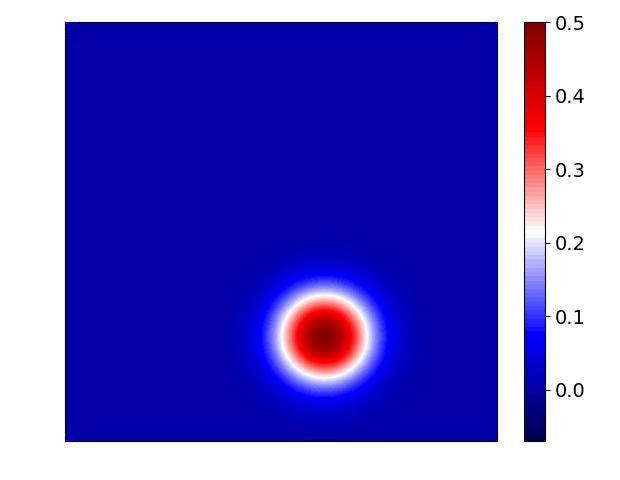}
         \subcaption{The target conductivity perturbation.}
     \end{subfigure}
     \hspace{1.5cm}
     \begin{subfigure}[t]{0.4\textwidth}
         \centering
         \includegraphics[width=\textwidth]{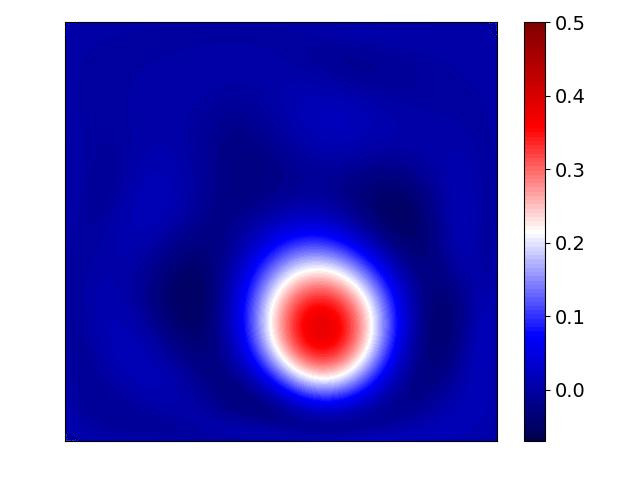}
         \subcaption{Reconstruction.}
     \end{subfigure}
     \hfill
     \caption{An SVD-based reconstruction of a smooth conductivity perturbation in a square for $M=16$, $\omega=1$ and $1\%$ of additive noise, resulting in $p=163$.}
     \label{fig:reco_square}
\end{figure}

\begin{figure}[htb]
     \centering
     \begin{subfigure}[t]{0.4\textwidth}
         \centering
         \includegraphics[width=\textwidth]{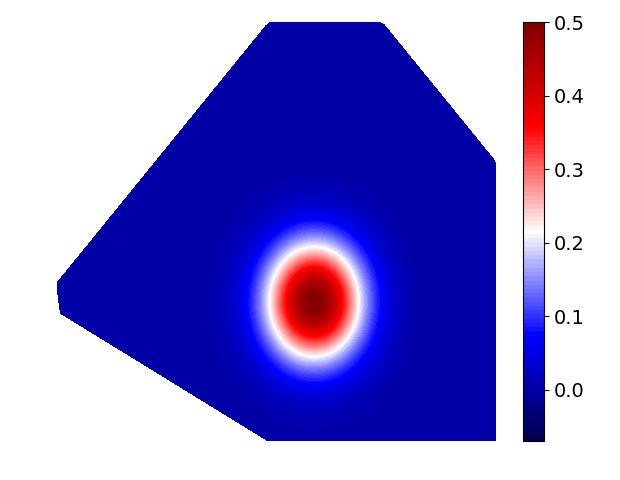}
         \subcaption{The target conductivity perturbation.}
     \end{subfigure}
     \hspace{1.5cm}
     \begin{subfigure}[t]{0.4\textwidth}
         \centering
         \includegraphics[width=\textwidth]{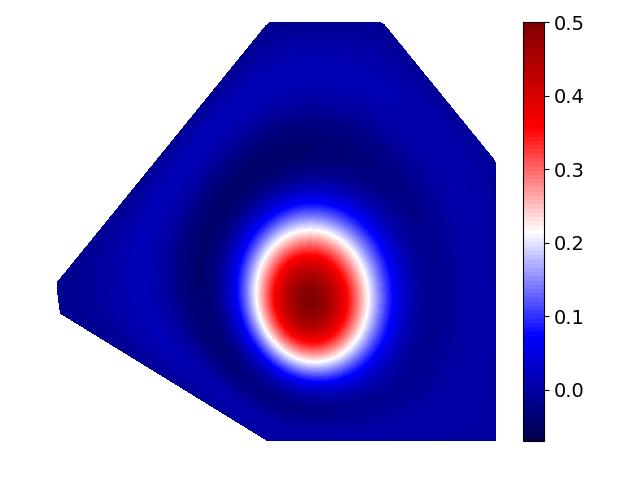}
         \subcaption{Reconstruction.}
     \end{subfigure}
     \hfill
     \caption{An SVD-based reconstruction of a smooth conductivity perturbation in a polygon for $M=16$, $\omega=1$ and $1\%$ of additive noise, resulting in $p=145$.}
     \label{fig:reco_hexagon}
\end{figure}

We then test the truncated SVD algorithm on the two polygonal domains shown in Figs.~\ref{fig:reco_square} and~\ref{fig:reco_hexagon}. The reason for only considering convex polygons as models for the imaged object $\Omega$ is the following: The derivative of the conformal mapping $\Psi'|_{\partial \Omega}$, needed in \eqref{eq:conformalfm} for defining the appropriate current basis on $\partial \Omega$, has a zero at each corner of $\partial \Omega$ with interior angle less than $\pi$ and a singularity at each corner with interior angle wider than $\pi$; see,~e.g.,~\cite[Section~3.2]{Trefethen80} for more details. Hence, by settling for convex polygons, one can simulate the required measurements on $\Omega$ without having to solve Neumann problems for singular boundary current densities. 

The considered conductivity perturbations are smooth Gaussian humps centered so far from the boundary that their supports are compactly contained in $\Omega$ within the numerical precision. The SVD-based algorithm, with $\omega = 1$, is applied to simulated data with 1\% of additive noise to produce the reconstructions on the right in Figs.~\ref{fig:reco_square} and~\ref{fig:reco_hexagon}. Although the algorithm finds the positions of the Gaussian humps, the reconstructions are not as localized as the target perturbations, and they also exhibit wave-like artifacts in the background. 

\subsection{CEM data for $D$} \label{sec:CEM_D}

\begin{figure}[htb]
     \centering
     \begin{subfigure}[t]{0.4\textwidth}
         \centering
         \includegraphics[width=\textwidth]{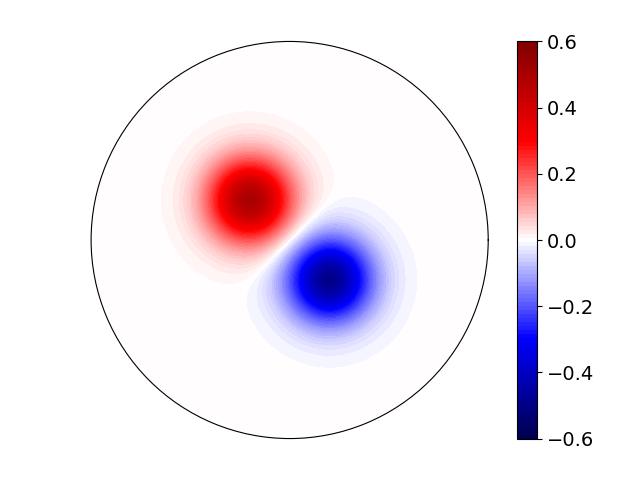}
         \subcaption{The target conductivity perturbation.}
         \label{fig:target_cem}
     \end{subfigure}
     \hspace{1cm}
     \begin{subfigure}[t]{0.4\textwidth}
         \centering
         \includegraphics[width=\textwidth]{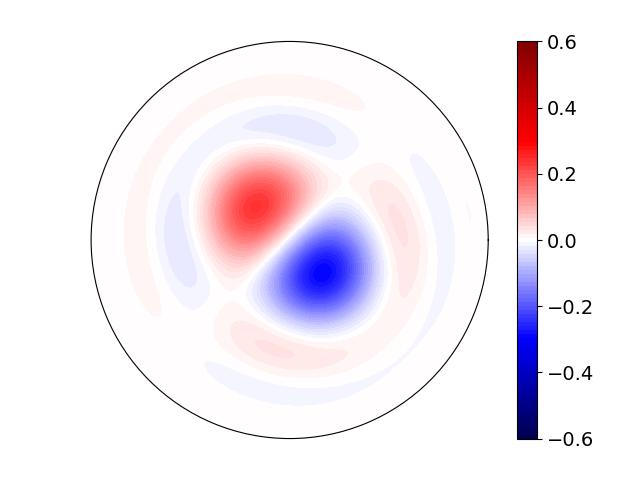}
         \subcaption{Reconstruction.}
         \label{fig:reco_cem}
     \end{subfigure}
     \hfill
     \caption{An SVD-based reconstruction of a smooth conductivity perturbation from simulated CEM data with 32 electrodes, $M = 16$, $\omega = 1$ and $1\%$ of additive noise, resulting in $p=108$.}
\end{figure}

Consider an electrode configuration with $32$ electrodes that are attached equiangularly to the unit disk and cover, in total, half of its boundary.\footnote{Although our model for the CEM in Section~\ref{sec:CEM} assumes an odd number of electrodes, the generalization of the approximation in \cite{Garde21} for an even number of electrodes is straightforward.} The contact impedance is set to $z =1$ on all electrodes. CEM measurements are simulated for $31$ linearly independent current patterns by a FEM with $33,281$ piecewise linear basis functions and appropriate refinements of the mesh at the electrodes (cf.~\cite{Vauhkonen97}); this is done for both the mere background conductivity and the background conductivity contaminated by the perturbation composed of two Gaussian humps with opposite signs depicted in Fig.~\ref{fig:target_cem}. As briefly explained in Section~\ref{sec:CEM}, the `measured' relative electrode potentials can then be used to construct an approximation for the data matrix $\{ a_{m,n}^{\rm NL}(\eta) \}_{m,n \in \Z_M'}$. See \cite{Garde21} for the details.

We choose $M=16$ as the maximal absolute value of the considered spatial Fourier frequencies and add $1\%$ of noise to the data. Observe that in addition to the synthetic measurement noise, there are two other significant sources of error in the data: the approximation between the CM and the CEM quantified by~\eqref{eq:CEM_approx} and the employment of `nonlinear' data for a linear reconstruction method. We employ the SVD-based algorithm with $\omega = 1$ to form the reconstruction presented in  Fig.~\ref{fig:reco_cem}. The reconstruction captures the basic characteristics of the target perturbation reasonably well. 

\subsection{Water tank measurements}

In the final experiment, we consider data from a right cylindrical water tank with 16 electrodes documented in~\cite{Hauptmann17}. Note that all measurements in \cite{Hauptmann17} are performed with electrodes that are homogeneous in the vertical direction and extend from the bottom of the tank exactly up to the water surface. The embedded inclusions are also homogeneous in the vertical direction and break the water surface, and they are either perfect insulators (plastic) or can be modeled as ideal conductors (metal). Due to the insulating,~i.e.,~homogeneous Neumann, boundary conditions at the bottom and top of the water layer, the measurements can thus be modeled by the two-dimensional CEM in the unit disk. Moreover, conductivity perturbations reconstructed by resorting to such a model could be presented in the SI units by appropriately scaling with the dimensions of the tank and the conductivity of the saline filling the tank. However, since it is not clear what the precise value for the conductivity of the saline is in \cite{Hauptmann17}, the presented reconstructions are to be understood as changes to the unit conductivity inside the unit disk without a proper interpretation in SI units.

\begin{figure}[htb]
     \centering
     \begin{subfigure}[t]{0.24\textwidth}
         \centering
         \includegraphics[width=\textwidth]{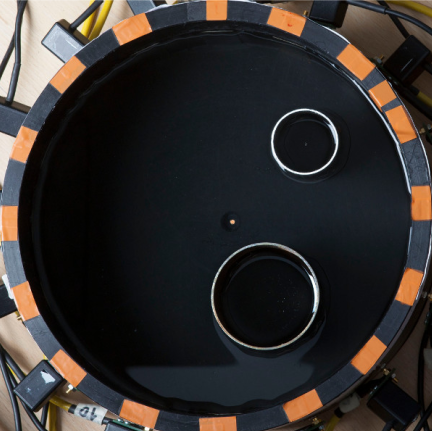}
     \end{subfigure}
\qquad \quad \ \
     \begin{subfigure}[t]{0.24\textwidth}
         \centering
         \includegraphics[width=\textwidth]{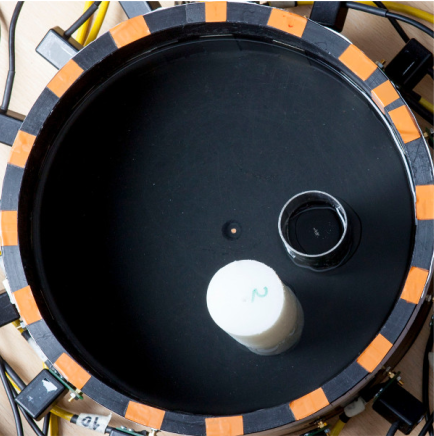}
     \end{subfigure}
     \qquad \quad \ \
     \begin{subfigure}[t]{0.24\textwidth}
         \centering
         \includegraphics[width=\textwidth]{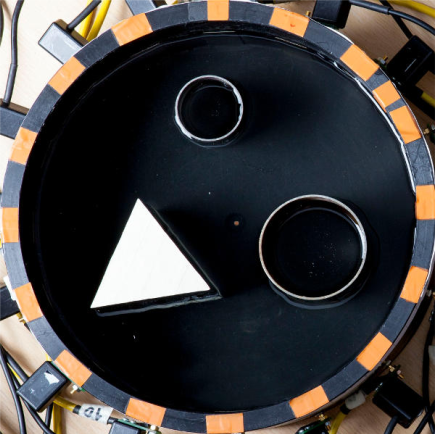}
     \end{subfigure}
     \quad \  \mbox{ }
     \\[7mm]
     \begin{subfigure}[t]{0.32\textwidth}
         \centering
         \includegraphics[width=\textwidth]{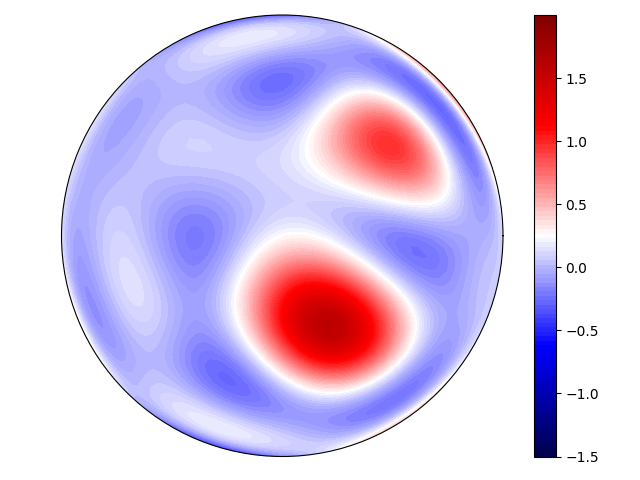}
         %\subcaption{Case 2.3}
     \end{subfigure}
     \hfill
     \begin{subfigure}[t]{0.32\textwidth}
         \centering
         \includegraphics[width=\textwidth]{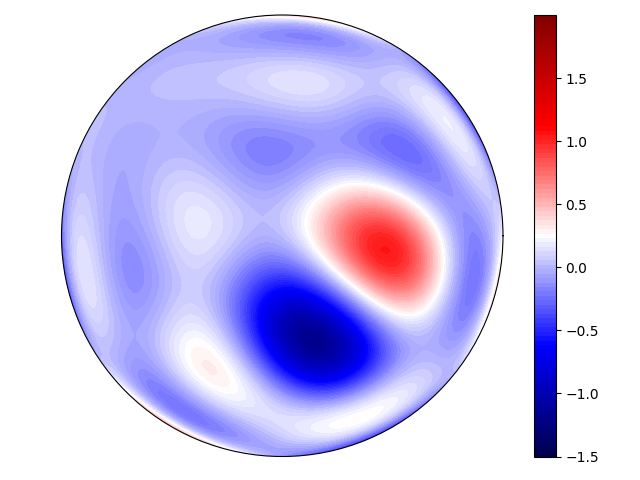}
         %\subcaption{Case 4.4}
     \end{subfigure}
     \hfill
     \begin{subfigure}[t]{0.32\textwidth}
         \centering
         \includegraphics[width=\textwidth]{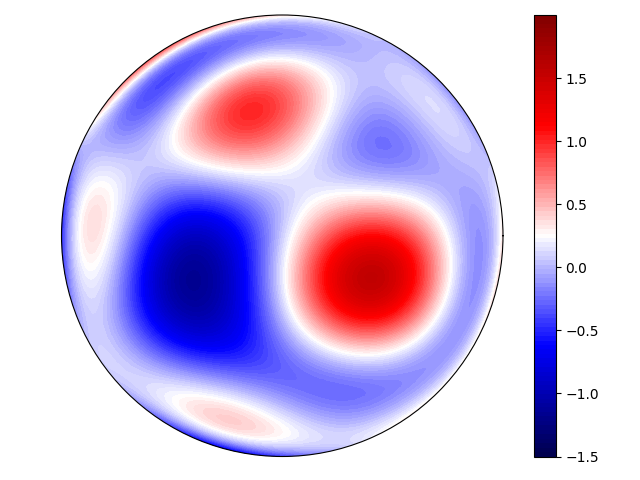}
         %\subcaption{Case 5.2}
     \end{subfigure}
     \hfill
     \caption{Three SVD-based reconstructions with $M=8$ from real-world water tank difference data measured at 16 equiangular electrodes. For each target, the `optimal' truncation index $p$ has been chosen according to a visual inspection. Top: the targets. Bottom: the reconstructions. }
     \label{fig:water_tank_recos}
\end{figure}

We employ the current pattern set 5 for the target cases~2.3, 4.4 and 5.2 as numbered in \cite{Hauptmann17}. The three targets are shown at top in Fig.~\ref{fig:water_tank_recos}, and the set 5, dubbed ``all against 1'' in \cite{Hauptmann17}, corresponds to feeding current through the right-most electrode and letting it out in turns through the other ones, resulting in 15 linearly independent current patterns. Instead of the set 5, one could as well use any other current patterns in~\cite{Hauptmann17}, as long as they form a basis for the space of mean-free vectors of length 16. Note that the measurements for the empty tank,~i.e.,~the case 1.0 in~\cite{Hauptmann17}, are also needed for building an approximation for the data matrix $\{ a_{m,n}^{\rm NL}(\eta) \}_{m,n \in \Z_M'}$, as reviewed in Section~\ref{sec:CEM_D}. In particular, the applicability of the approximation \eqref{eq:CEM_approx} is based on an assumption that the contacts at the electrodes are (approximately) the same for all considered target cases.

We choose $M=8$ and apply the SVD-based reconstruction algorithm to the three sets of data corresponding to the configurations on the top row of Fig.~\ref{fig:water_tank_recos}. Since we do not know the noise level (or the noise model) exactly, we decide to tune the truncation index $p$ via trial and error so that each of the three reconstructions is as good as possible according to a subjective visual inspection. The resulting reconstructions are shown on the bottom row of Fig.~\ref{fig:water_tank_recos}. According to our noise model and the Morozov discrepancy principle, the used regularization parameters correspond to noise levels between $0.1\%$ and $1\%$. Although the reconstructions in Fig.~\ref{fig:water_tank_recos} are not based on any systematic way of choosing the level of regularization, they in any case demonstrate that the presented algorithm can under optimal circumstances produce reasonable reconstructions in practical EIT.

\section{Concluding remarks}\label{sec:conclusion}

This paper introduced, implemented and numerically tested two regularized variants of the direct reconstruction method for the linearized CM of two-dimensional EIT presented in~\cite{Garde22}. Both algorithms exploit the possibility to present the linearized problem in the unit disk with the help of a block diagonal (infinite) matrix with lower triangular blocks, which enables tackling the triangular subsystems separately in the inversion process. More precisely, one of the regularized algorithms is based on a semi-heuristic truncation of the triangular subsystems and the other one on truncating their SVDs. According to our numerical tests, the latter variant is more stable and capable of producing high-quality reconstructions from accurate simulated data, as well as reasonably good reconstructions from noisy measurements on convex polygonal domains and even from real-world water tank data. In particular, ignoring the linearization error does not seem to significantly impede the practical applicability of the introduced reconstruction method.

\subsection*{Acknowledgments}

This work is supported by the Academy of Finland (decisions 353080, 353081, 359181) and the Aalto Science Institute (AScI). HG is supported by grant 10.46540/3120-00003B from Independent Research Fund
Denmark \textbar\ Natural Sciences.

\bibliographystyle{plain}
\bibliography{numlincal-refs}

\end{document}